\documentclass[MSNbibl,aap,citesort,dvips]{arximspdf}
\usepackage{mathbh}
\usepackage{graphicx}

\doi{10.1214/10-AAP739}
\volume{21}
\issue{4}
\pubyear{2011}
\firstpage{1537}
\lastpage{1567}

\makeatletter

\newcommand{\fraca}[2]{{(#1/#2)}}
\newcommand{\LL}{\ll}
\newcommand{\binom}[2]{\pmatrix{#1\cr#2}}

\newcommand{\mb}{\mathbf}
\newcommand{\eps}{\varepsilon}
\newcommand{\w}{W_1}
\newcommand{\ov}{\overline}
\renewcommand{\overline}{\bar}
\newcommand{\Lip}{\operatorname{Lip}_1}
\newcommand{\Lipd}{\operatorname{Lip}^{\Delta}_1}
\newcommand{\cal}{\mathcal}
\newproclaim{example}{Example}
\newproclaim{remark}{Remark}
\newtheorem{theorem}{Theorem}
\newtheorem{proposition}{Proposition}
\newtheorem{corollary}{Corollary}
\newtheorem{lemma}{Lemma}
\makeatother

\begin{document}
\begin{frontmatter}

\title{Scaling limits for continuous opinion dynamics~systems}
\runtitle{Scaling limits for continuous opinion dynamics systems}

\begin{aug}
\author[A]{\fnms{Giacomo} \snm{Como}\corref{}\ead[label=e1]{giacomo@mit.edu}}
\and
\author[B]{\fnms{Fabio} \snm{Fagnani}\ead[label=e2]{fabio.fagnani@polito.it}}
\runauthor{G.~Como and F.~Fagnani}
\affiliation{Massachusetts Institute of Technology and Politecnico di Torino}
\address[A]{Laboratory for Information\\ \quad and Decision Systems\\
Massachusetts Institute of Technology\\ 77 Massachusetts Ave\\
Cambridge, Massachusetts 02139\\ USA\\
\printead{e1}} 
\address[B]{Dipartimento di Matematica\\ Politecnico di Torino\\ Corso
Stati Uniti 24\\
10129 Torino\\
Italy\\
\printead{e2}}
\end{aug}

\received{\smonth{3} \syear{2010}}
\revised{\smonth{9} \syear{2010}}

%
\begin{abstract}
Scaling limits are analyzed for stochastic continuous opinion dynamics systems,
also known as gossip models. In such models, agents update their
vector-valued opinion to a convex combination (possibly agent- and
opinion-dependent) of their current value and that of another observed
agent. It is shown that, in the limit of large agent population
size, the empirical opinion density concentrates, at an
exponential probability rate, around the solution of a
probability-measure-valued ordinary differential equation describing
the system's mean-field dynamics.
Properties of the associated initial value problem are studied. The asymptotic
behavior of the solution is analyzed for bounded-confidence opinion dynamics,
and in the presence of an heterogeneous influential environment.
\end{abstract}

%

\end{frontmatter}

\section{Introduction}\label{sectintro}
In this paper, we undertake a rigorous mathematical analysis of a
family of stochastic dynamical systems proposed as opinion
dynamics models in the recent literature: see, for example,
\cite{CastellanoFortunatoLoreto}, Section~III,~\cite{Lorenz}, and
references therein. Here, we shall focus on the so-called ``gossip''
models, where the information propagation, as the name suggests,
takes place through pairwise interactions. These models have been
proposed in other scientific areas, for instance, as aggregation
and estimation algorithms in sensor and robotic networks (see,
e.g., \cite{BoydGossip,Shah}), or as models for aggregation and
clustering in biological systems (see, e.g., \cite{Edelstein}).

One of the simplest gossip model can be described as follows. Each
agent~$a$ of a population ${\cal A}$ of finite size $n:=|{\cal
A}|$ possesses an initial belief/opinion modeled as a vector
$X^{a}_0\in{\mathbb{R}}^d$. Agents are activated according to independent
Poisson processes in continuous time.\footnote{Analogous versions of
this model have been presented in the literature with agents'
activations occurring in discrete time.} If agent $a$ is activated at
time $t$, her opinion jumps from its current value $X^a_{t^-}$ to
a new value $X^{a}_{t}=\ov\omega X^{a}_{t^-}+\omega
X^{b}_{t^-}$ where $b$ is another agent sampled from
$\mathcal A$, and $\omega=1-\ov\omega\in[0,1]$ is a parameter modeling
how much
agent $a$ trusts the opinion of agent $b$. In general, the
conditional distribution of $b$ over the agent population may
depend on the activated agent $a$ (the support of such
distribution representing the out-neighborhood of $a$ in an
underlying ``social network'' structure), while the parameter
$\omega$ may depend on the interacting agents, $a$ and $b$, as
well as on their current opinions, $X_{t^-}^a$ and $X_{t^-}^b$.

Fundamental theoretical issues concern the behavior of such models for
large $t$ and large $n$.
Rather then in the single opinions' behavior, one is interested in the
emerging collective behavior of the population.
Typical questions include whether a consensus is eventually achieved or
rather disagreement persists, and, more in general, whether an
asymptotic distribution of opinions exists, what it looks like, and how
long it takes the system to approach~it.\looseness=-1

The simplest case is when the Poisson processes are all of unitary
rate, the conditional distribution of the
observed agent is uniform over $\mathcal A$ whichever agent is
activated, and the parameter $\omega$ is fixed and the same for
all agents, independently of their current opinions. In this case,
the model is linear and can be studied in full detail: it corresponds
to the asymmetric gossip model in \cite{FagnaniZampieri}. The
basic fact is that (if $\omega\in\,]0,1[$), almost surely, all
$X^{a}_{t}$ converge, as $t\to+\infty$ (and for any fixed $n$),
to a consensus random value $\xi$ which has expected value
${\mathbb{E}}(\xi)=n^{-1}\sum_{a}X^{a}_0$. Convergence is
exponentially fast
\cite{FagnaniZampieri2}:
\[
 \mathbb{E} \biggl [n^{-1}\sum_{a}|X_t^{a}-\xi|^2 \biggr]\leq
2n^{-1}\sum
_{a}|X_0^{a}|^2\exp(-Ct) ,
\]
where $C=-n\ln(1-2n^{-1}\omega\bar\omega
-2n^{-2}\omega^2)$. The variance of $\xi$ can be estimated as
\[
\operatorname{Var}[\xi]\leq\frac{\omega}{\omega+\bar\omega
n}n^{-1}\sum
_{a}|X_0^{a}|^2 .
\]
Moreover, using the techniques in \cite{FagnaniZampieri}, one can
easily prove a concentration result of type
\[
{\mathbb{P}}\bigl(|X_t^a-{\mathbb{E}}(X_t^a)|\geq\varepsilon\bigr)\leq
\exp
(-K\varepsilon^2n/t ) .
\]
Essentially, this shows that, as $n$ grows large,
and $t/n$ tends to $0$, each agent's opinion $X^a_t$ concentrates
around a deterministic dynamics
converging to ${\mathbb{E}}(\xi)$ as $\exp(-2\omega\bar\omega t)$.
It is this type of results which we would like to extend to more
general models.

A particularly interesting setting is the homogeneous-population,
state-depend\-ent model, that is, when the parameter $\omega$ is
independent of the identity of the interacting agents, but does depend
on their
current opinions. The case
%
\begin{equation}\label{boundedconfidence}\omega
=\omega(X^a_{t^-},X^b_{t^-})=
\cases{\displaystyle\omega_0,&\quad  if
 $|X^{a}_{t^-}-X^{b}_{t^-}|\le R$,\cr\displaystyle
 0,&\quad  if
  $|X^{a}_{t^-}-X^{b}_{t^-}|> R$,
}
\end{equation}
where $R>0$, and $\omega_0\in\,]0,1[$, is known as the
Deffuant--Weisbuch model \cite{Krause,Deffuant,Lorenz0} of bounded
confidence opinion dynamics: agents with opinions too far apart do not
trust each other, hence they do not interact.
Another case is the so-called Gaussian interaction kernel
%
\begin{equation}\label{Gaussiandecay}
\omega=\omega(X^a_{t^-},X^b_{t^-})
=\omega_0\exp(-|X^{a}_{t^-}-X^{b}_{t^-}|^2/\sigma^2 )
,
\end{equation}
a similar form of which was considered in \cite{Deffuant2}.
Observe that, in these models, the dynamics of the network and of the
opinions become intertwined.
In fact, these models are nonlinear and, to the best of the authors'
knowledge, the
only theoretical result \cite{Lorenz0} is that, if $\omega\in
\{0\}\cup[\omega_0,1]$ for some $\omega_0>0$,
each $X^{a}_{t}$ converges, as $t$ grows large, to a limit random value
$\xi^{a}$. Numerical
simulations show the asymptotic emergence of opinion clusters whose
number and structure depends on the initial condition but seems to be
stable for large $n$. However, there is no theoretical result regarding
concentration and scaling limits for any state-dependent model.

On the other hand, for the case when the parameter $\omega$ depends on
the agents, as well as on their opinions, no theoretical result is
available in the literature. Some of these heterogenous models have
been considered in \cite{Weisbuch,Fortunato,Lorenz1,Lorenz2} where,
though, only numerical simulations have been presented. Such
heterogeneous population models are going to play
a very important role in opinion dynamics because they are the
natural model to represent more realistic populations with agents
having different attitude to change opinion, and interacting only with
agents in their social neighborhood.

In this paper, we study general state-dependent gossip models for
large~$n$. We shall consider both the case of a homogeneous population, and
of a~heterogeneous one consisting of two classes of agents: ``standard''
agents, which keep on updating their opinions as a result of
interactions with the whole population, and ``stubborn'' agents whose
opinions are never updated~\cite{disagreement}. The latter case can be
modeled as a homogeneous population model with an exogenous input
describing the influence of the stubborn agents' opinions on the
standard agents' ones, and interpreted as~a, typically heterogeneous,
``influential environment.'' We believe that many more general
heterogeneous models can be studied with our approach. This will be
done in a forthcoming paper where also models with interactions of
nongossip type will be considered.

In our analysis, we shall adopt an Eulerian viewpoint: instead of
studying the evolution of the single agents' opinions, we shall neglect
the agents' identities, and study the dynamics of the corresponding
empirical opinion densities. We shall argue that the deterministic
mean-field dynamics obtained in the limit of large $n$ is governed by
an ordinary differential equation (ODE) on the space of probability
measures over the opinion set, presented in Section~\ref{sectEuler}.
As proved in Section~\ref{sectsolutions}, the initial value problem
associated to the mean-field dynamics always admits a unique global
solution. Moreover, at any finite time, its solution is absolutely
continuous with respect to Lebesgue's measure, provided that so does
the initial condition, and that some mild technical conditions are
satisfied by the interaction kernel.

The asymptotic behavior in time of the mean-field dynamics is analyzed
in Section~\ref{sectexamples} for the state-independent heterogeneous
case, and for the generally state-dependent homogeneous case. In both
cases, we prove weak convergence to an equilibrium distribution, which
typically does not consist of a single Dirac's delta. For the
state-independent heterogeneous model, we show that the equilibrium
opinion distribution is independent of the initial condition, and is
uniquely characterized by its moments, which can be computed by
recursively solving a lower-triangular infinite linear system. On the
other hand, we prove that the equilibrium opinion distribution in the
bounded-confidence model is a convex combination of Dirac's deltas.
Such deltas represent opinion clusters, and their number and position
depend on the initial condition.\footnote{Proofs of similar results
showing convergence of various variants of the bounded confidence
opinion dynamics to opinion clusters have appeared in \cite
{Lorenz0,Lorenz1,CanutoFagnaniTilli,BlondelHendrikxTsitsiklis1,BlondelHendrikxTsitsiklis2}.}
Our results provide fundamental insight into two basic mechanisms which
have been proposed by social scientists in order to explain persistent
disagreement in the society \cite{Axelrod}, namely heterogeneity of
the social environment, and homophily leading to global fragmentation.

Finally, in Section~\ref{sectconcentration}, we prove that the
finite-population stochastic system concentrates around the
deterministic mean-field dynamics, as the population size grows, at an
exponential probability rate. We apply here a martingale argument (see,
e.g., \cite{Wormald} for the finite-dimensional case) and obtain a~%
result in the Kantorovich--Wasserstein metric \cite
{AmbrosioGigliSavare,Villani}. The technical assumption in our results
is that the, possibly stochastic, dependence of the weight~$\omega$ on
the opinions is Lipschitz-continuous. Hence, the case (\ref
{boundedconfidence}) is not covered by our theory. This is not a
relevant drawback since one can consider suitable Lipschitz
approximations of (\ref{boundedconfidence}); on the other hand, we
believe that this is just a technical question and that the result
should remain valid for a~larger class of functions.

We conclude this section with a brief overview of some related work.
A~special instance of the measure-valued ODE analyzed in the present
paper has
already been proposed in \cite{BenNaim} for probability densities (in
this case it becomes an integro-differential equation), but
with no proof of either well-posedness or concentration of the
stochastic finite system.
In \cite
{BlondelHendrikxTsitsiklis1,CanutoFagnaniTilli,BlondelHendrikxTsitsiklis2},
deterministic, bounded-confidence, opinion dynamics models with
possibly a continuum of agents have been studied both in discrete
and continuous time. In particular, the continuous-time opinion
dynamics studied in \cite{CanutoFagnaniTilli} is governed by a partial
differential equation in the space of probability measures, while the
work~\cite{BlondelHendrikxTsitsiklis2} deals with the equivalent
dynamics, in dimension one, of the cumulative distribution functions.
In both works, the agents' opinions have continuous trajectories, and
the corresponding generator of the opinion density dynamics is local.
In contrast, in the model analyzed in the present paper, the opinion
trajectories are discontinuous (in fact, piece-wise constant), and the
induced mean-field dynamics is driven by a nonlocal operator. As shown
in Section~\ref{sectinfluential}, the bounded-confidence mean-field
dynamics studied here has a~qualitatively similar behavior to the
solution of the partial differential equation of~\cite
{CanutoFagnaniTilli,BlondelHendrikxTsitsiklis2}. It is also worth
mentioning the work \cite{HaLiu}, where mean-field limits have been
analyzed for the flocking dynamics of Cucker and Smale \cite
{CuckerSmale1,CuckerSmale2}. Finally, only at the end of their work the
authors have become aware that an approach very similar to the one in
this paper has been undertaken in \cite{Remenik}, based on results in
\cite{Fournier}.

\section{Problem setting and main results}\label{sectprobsettingmainresults}
In this section, we formally state the model and present our main results.

Before proceeding, let us establish some notation to be followed
throughout the paper.
For $x,y\in{\mathbb{R}}^d$, for some $d\in{\mathbb{N}}$, $|x-y|$
and $x\cdot y$ will
denote their Euclidean distance, and scalar product, respectively.
The indicator function of a set $A$ will be denoted by $\mathbh{1}_A$,
that is,
$\mathbh{1}_A(x)=1$ if $x\in A$, and $\mathbh{1}_A(x)=0$ if $x\notin A$.
Given an open subset $\mathcal X\subseteq{\mathbb{R}}^d$, we denote
by $\mathcal B(\mathcal
X)$ its Borel $\sigma$-algebra, and by $\mathcal M( \mathcal X)$ the space
of finite signed Borel measures on $\mathcal X$, equipped with the topology
of weak-$*$ convergence,
while $\mathcal M^+(\mathcal X)\subseteq\mathcal M(\mathcal X)$
denotes the closed convex
cone of Borel
nonnegative measures, and $\mathcal P(\mathcal X)\subseteq\mathcal
M^+(\mathcal X)$
the simplex of probability measures over $\mathcal X$. The space of real-valued
continuous bounded (resp., compact-supported, vanishing at
infinity) functions on~$\mathcal X$, equipped with the supremum norm
$\|\varphi\|_{\infty}:=\sup\{|\varphi(x)|\dvtx  x\in\mathcal X \}
$, will
be denoted by $\mathcal C_b(\mathcal X)$ [resp., $\mathcal C_c(\mathcal
X)$, $\mathcal
C_0(\mathcal X)$].
The Dirac delta measure centered in $x\in\mathcal X$ will be denoted
by $\delta_x$. For $\mu\in\mathcal M(\mathcal X)$  and $\varphi\in
\mathcal
C_b(\mathcal X)$, we shall write $\langle\mu,\varphi\rangle$ for the
integral $\int\varphi(x)\,d\mu(x)$, with the convention that,
whenever not explicitly indicated, the domain of integration is
assumed to be the entire space $\mathcal X$. The total variation of
$\mu\in\mathcal M(\mathcal X)$ will be denoted by $\|\mu\|$.
The symbol $\lambda$ will denote Lebesgue's measure on $\mathcal X$,
$\mu\LL\lambda$ will stand for absolute continuity, and
$ d\mu/ d\lambda$ for the Radon--Nikodym derivative, of $\mu$
with respect to $\lambda$. Finally, we shall denote by $\mathcal
P_1(\mathcal
X):=\{\mu\in\mathcal
P(X)\dvtx  \int|x|\,d\mu(x)<+\infty\}$ the metric space of probability
measures with finite first moment, equipped with the order-$1$
Kantorovich--Wasserstein distance. The latter is defined by
$\w(\mu,\nu):=\inf\{\int\!\!\int|x-y|\,d\xi(x,y) \}$,
where the
infimization (which is in fact a minimization \cite
{AmbrosioGigliSavare,Villani}) runs over all couplings of $\mu$ and
$\nu$,
that is, joint probability measures $\xi\in\mathcal P(\mathcal
X\times\mathcal X)$
having marginals given by $\mu$, and $\nu$, respectively.

\subsection{Stochastic models of continuous opinion dynamics} \label
{sectmodel}
The present paper is concerned with continuous opinion dynamics systems.
Agents belong to a finite population $\mathcal A$ of cardinality
$|\mathcal A|=n$.
At time $t\in{\mathbb{R}}^+$ each agent\vadjust{\eject} $a\in\mathcal A$ maintains
an opinion
$X^{a}_t\in\mathcal X$, where $\mathcal X\subseteq{\mathbb{R}}^d$ is
an open set.
The vector of the opinions will be denoted by $X_t:=\{X^{a}_t\dvtx  a\in
\mathcal A\}\in\mathcal X^{\mathcal A}$.

We shall assume the initial opinions $X_0$ to be a collection of
independent and identically distributed random variables, the law
of each $X^{a}_0$ given by some $\mu_0\in\mathcal P({\mathbb
{R}}^d)$. The
trajectories of the opinion profile vector $\{X_t\dvtx  t\in{\mathbb
{R}}^+\}$ are
right-continuous and evolve according to the following jump Markov
process: Agents have clocks which tick at the times of independent
rate-$1$ Poisson processes. If her clock ticks at time $t$, agent
$a$ updates her opinion~$X^{a}_{t^-}$ to a new value $X^{a}_{t}$
which depends on the observation of the current opinion of some
other agent and of her own one. In particular, she observes the opinion of
some other agent $b$ sampled uniformly from $\mathcal A$, and then
updates her opinion to a random value $X^{a}_t$, which has
conditional probability law
$\kappa( \cdot|X^{a}_{t^-},X^{b}_{t^-})$. Here
$\kappa( \cdot| \cdot, \cdot)$ is a stochastic kernel,
that is, for all $x,y\in\mathcal X$, $\kappa( \cdot|x,y)$ is a
probability measure on $\mathcal X$, and $(x,y)\mapsto\kappa(B|x,y)$ is
a measurable map from $\mathcal X\times\mathcal X$ to $[0,1]$, for all
measurable sets $B\subseteq\mathcal X$. We shall refer to $\kappa$ as the
interaction kernel of the model.
We shall assume that the above stochastic process is defined on some
filtrated probability space
$(\Omega,\{\mathcal F_t\}_{t\in{\mathbb{R}}^+},{\mathbb{P}})$, and
denote by $0=T_0<
T_1<T_2<\cdots,$ the times at which some opinion update occurs (strict
inequalities holding almost surely).
Observe that $\{T_{k+1}-T_k\dvtx  k\in{\mathbb{Z}}^+\}$ is a family of
independent rate-$n$ Poisson random variables.

In most of the models considered in the literature, the interaction
kernel is a convex combination of type:
$\kappa( \cdot|x,y)=\alpha\kappa^i( \cdot|x,y)+\ov\alpha
\kappa^e( \cdot|x)$ where $\alpha=1-\ov\alpha\in[0,1]$ and
where $\kappa^i( \cdot|x,y)$ is a probability measure
concentrated on the interval connecting $x$ and $y$ while
$\kappa^e( \cdot|x)$ is a probability measure concentrated on
the segment connecting $x$ to some random point $z$. More
specifically, $\mathcal X\subseteq{\mathbb{R}}^d$ is a convex open
set containing
the support of the initial condition, and there exists two scalar
stochastic kernels $\theta^i( \cdot| \cdot, \cdot)$ and
$\theta^e( \cdot| \cdot, \cdot)$ from $\mathcal X\times
\mathcal X$
to $[0,1]$ such that
\[
\kappa^i( \ov\omega x+\omega
y |x,y)=\theta^i(\omega|x,y) , \qquad  \kappa^e( \ov\upsilon
x+\upsilon
z |x)=\int \theta^e(\omega|x,z)\,d\psi(z),
\]
where
$\ov\omega=1-\omega$, $\ov\upsilon=1-\upsilon$  and $\psi\in
\mathcal P(\mathcal X)$.
This models a situation in which, with probability $\alpha$, the
activated agent updates her opinion towards a convex combination
of her current opinion $x$ and the opinion $y$ of an observed
agent. The weight $\omega$ in such a convex combination measures
the confidence that the activated agent has on the observed
opinion of another agent, and is assumed to depend, through the
stochastic kernel $\theta^i( \cdot| \cdot, \cdot)$, on
both the activated and the observed agent's opinions, $x$ and $y$.
On the other hand, with probability $\ov\alpha$, the activated
agent observes an external signal $z$, sampled from a probability
distribution $\psi$, playing the role of an exogenous source of
influence, or influential environment, and she updates her opinion
toward a convex combination of her current opinion $x$ and the
observed signal $z$. The dependence of the weight $\upsilon$ of such
convex combination is captured by the stochastic kernel
$\theta^e( \cdot| \cdot,\cdot)$. A useful equivalent way
to characterize the interaction kernel $k$ described above is
through its action on continuous test functions:
%
\begin{eqnarray}\label{interacting
kernels}
\langle\kappa( \cdot|x,y),\varphi\rangle&=&
\alpha\int \varphi(\ov\omega x+\omega
y )\,d\theta^i(\omega|x,y)\nonumber
\\[-8pt]
\\[-8pt]
&&{}+
\ov\alpha\int\!\!\int \varphi(\ov\upsilon x+\upsilon
z)\,d\theta^e(\upsilon|x,z)\,d\psi(z)
\nonumber
\end{eqnarray}
for all $\varphi\in\mathcal
C_0(\mathcal X)$.

\begin{example}[(Gossip model with heterogeneous influential
environment)]\label{exampleexternalgossip}
Assume that the stochastic kernel $\kappa( \cdot| \cdot,
\cdot)$ has the form (\ref{interacting kernels}), with constant weights:
\[
\theta^i( \cdot|x,y)=\delta_{\omega}( \cdot) , \qquad  \theta
^e( \cdot|x,y)=\delta_{\upsilon}( \cdot)
\]
for some fixed confidence weights $\omega,\upsilon\in[0,1]$. This models
a homogeneous population whose opinion dynamics alternates internal
gossip updates to interactions with a static, external influential
environment. Internal gossip steps occur with probability $\alpha$,
and involve a uniformly sampled agent $a$ updating her opinion to a
convex combination, with trust parameter $\omega$, of her current
value and the one of another uniformly sampled agent $b$. Interactions
with the external environment occur with probability $\ov\alpha$, and
involve a uniformly sampled agent $a$ updating her opinion to a convex
combination, with trust parameter $\omega$, of her current value and
an external signal~$z$ sampled from a static distribution $\psi( dz
)$. This model has been analyzed in \cite{disagreement} for finite,
possibly inhomogeneous populations. The mean-field limit of this model,
with homogeneous population, will be analyzed in detail in Section~\ref
{sectinfluential}.
\end{example}

\begin{example}[(Bounded confidence opinion dynamics)]\label
{exampleboundedconfidence}
Consider the case when $\kappa( \cdot| \cdot, \cdot)$ is
in the form (\ref{interacting kernels}) with $\alpha=1$, and trust
parameter distribution $\theta^i( \cdot|x,y)$ supported on
$[0,\omega_0]$ for some $\omega_0\in[0,1[$. The case when $\theta
^i( \cdot|x,y)=\delta_{\omega(x,y)}$, where $\omega(x,y)$ is a
nonincreasing function of the distance $|x-y|$ can be consider to model
a homophily mechanism whereby agents are more likely to interact with
others which have similar opinions. In particular, the case when
$\omega(x,y)=0$ for all $|x-y|>R$, for some finite $R>0$, is usually
referred to as bounded confidence opinion dynamics~\cite
{Deffuant,Lorenz0}, and the minimum such~$R$ as the confidence
threshold. The special case $\omega(x,y)=\omega_0\mathbh{1}_{[0,R]}(|x-y|)$
corresponds to the Deffuant--Weisbuch model~\cite
{Deffuant,BenNaim,Lorenz}. The mean-field limit of the bounded
confidence opinion dynamics model will be analyzed in detail in
Section~\ref{sectboundedconfidence}.
\end{example}

While for most of the results of our paper we shall not need the
interaction kernel $\kappa$ to have the specific form
(\ref{interacting kernels}), we shall focus on kernels of this
form in Section~\ref{sectexamples} when proving asymptotic properties
of the
solution of the corresponding measure-valued ODE.

\begin{remark}
The models considered in the cited literature often assume the
interaction to be symmetric: when agent $a$ is activated and interacts
with agent $b$, both agents update their opinions. This symmetric
model may be more suitable in certain applicative contexts, the
asymmetric one in some others. However, while for finite
population sizes some of the properties of the two models differ
(e.g., in the symmetric model the average of the opinions
is preserved, while this is not necessarily the case for the
asymmetric model~\cite{FagnaniZampieri}), all the results and
proofs of this paper hold, with minor changes, for the
symmetric model too.
\end{remark}

\subsection{The Eulerian viewpoint and main results}\label{sectEuler}
As the main interest is in the global behavior of the opinion
dynamics system, rather than on that of the single agents'
opinions, it proves convenient to adopt an Eulerian viewpoint,
studying the evolution of the empirical densities of the agents'
opinions. Formally, this is accomplished by considering the random
flow of probability measures
\[
\mu^{n}_t:=\frac1n\sum_{a\in\mathcal A}\delta_{X^{a}_t}\in
\mathcal P(\mathcal
X) , \qquad  t\in{\mathbb{R}}^+ .
\]
This is a $\mathcal P(\mathcal X)$-valued process whose trajectories are
piecewise constant and right continuous.
In particular, one has
\[
\mu^n_t=M_{k}    \qquad\forall t\in[T_k,T_{k+1}[ ,\ k\in
{\mathbb{Z}}^+ ,
\]
where $\{M_{k}\dvtx  k\in{\mathbb{Z}}^+\}$ is a $\mathcal P(\mathcal
X)$-valued Markov chain.

In order to describe the dynamics of the $M_k$'s it is
useful to consider the operator $F\dvtx \mathcal M^+(\mathcal X)\to\mathcal
M^+(\mathcal
X)$, defined by
\[
F(\mu)(B):=\int\!\!\int\kappa(B|x,y)\,d\mu(x)\,d\mu(y)  \qquad \forall B\in\mathcal B(\mathcal X) .
\]
Equivalently, one can write
%
\begin{equation}\label{defop}
\langle
F(\mu),\varphi\rangle:=\int\!\!\int\!\!\int\varphi(z )\,d
\kappa(z|x,y)\,d\mu(x)\,d\mu(y)
\end{equation}
for all $\varphi\in\mathcal C_0(\mathcal X)$.
When $\mu$ is a probability measure, then $F(\mu)$ may be interpreted
as the conditional distribution of the new opinion formed as a~result
of the first interaction occurring, given that the current empirical
opinion density is $\mu$. In fact, the opinions $x$ and $y$ of two
agents $a$ and~$b$, randomly sampled, independently and uniformly, from
the agent population, have conditional joint distribution $ d\mu
(x)\,d\mu(y)$, and hence the new opinion $z$ formed as a result of
their interaction has conditional distribution $ d\kappa(z|x,y)\,d
\mu(x)\,d\mu(y)$.


It is immediate to verify that
\[
 \mathbb{E} [\langle M_{k+1},\varphi\rangle| M_{k}
]=(1-n^{-1} )\langle
M_{k},\varphi\rangle+n^{-1}\langle F(M_{k}),\varphi\rangle
\]
for all $\varphi\in\mathcal C_0(\mathcal X)$  and $k\in{\mathbb{Z}}_+$.
One may rewrite this in the form
%
\begin{equation}\label{sampleeq2} \langle M_{k+1},\varphi\rangle
-\langle
M_{k},\varphi\rangle= n^{-1}\langle F(M_{k})-M_{k},\varphi\rangle
+n^{-1}\langle\Lambda_{k+1},\varphi\rangle,
\end{equation}
where the random
signed measure $\Lambda_{k+1}$ satisfies
%
\begin{equation}\label{Deltaprop} \quad
{\mathbb{E}}[\Lambda_{k+1}|\mathcal
F_{T_k} ]=0 , \qquad  \|\Lambda_{k+1}\|\le
n\|M_{k+1}-M_k\|+\|F(M_k)-M_k\|\le4 .
\end{equation}
Equation
(\ref{Deltaprop}) implies that
$\{\langle\Lambda_{k},\varphi\rangle\dvtx  k\in{\mathbb{N}}\}$ is a
sequence of bounded
martingale differences, which can be thought as ``noise.'' This
suggests to think of the stochastic process $\{M_k\dvtx  k\in{\mathbb
{Z}}^+\}$ as
of a noisy
discretization, or Euler approximation in the numerical analysis
language, of the probability-measure-valued ODE
%
\begin{equation}\label{ODE}
\frac{ d}{ dt}\mu_t=F(\mu_t)-\mu_t
\end{equation}
with stepsize $1/n$.
We shall refer to a solution of (\ref{ODE}) as the mean-field dynamics
of the system.

More precisely, we shall define a solution of (\ref{ODE})
to be a family of probability measures $\{\mu_t\dvtx  t\in[0,+\infty)\}$
such that, for every function $\varphi\in\mathcal C_0(\mathcal
X)$, the real-valued map $t\mapsto\langle\mu_t,\varphi\rangle$ is
differentiable on ${\mathbb{R}}^+$, and satisfies
%
\begin{equation}\label{weaksol}\frac{ d}{ dt}\langle\mu
_t,\varphi\rangle=\langle
F(\mu_t),\varphi\rangle-\langle\mu_t,\varphi\rangle
\end{equation}
for every
$t>0$. The main result of this paper, stated below, guarantees
that~(\ref{ODE}) admits a unique solution $\{\mu_t\}$, and that
the stochastic process~$\{\mu^n_t\}$ concentrates around
$\{\mu_t\}$ exponentially fast in $n$.
\begin{theorem}\label{theomain}Let $\mu\in\mathcal P(\mathcal X)$ be
arbitrary. Then:
\begin{longlist}[(b)]
\item[(a)] There exists a unique
solution $\{\mu_t\dvtx  t\in{\mathbb{R}}^+\}$ of (\ref{ODE}) with initial
condition $\mu_0=\mu$;
%
\item[(b)] If $\mathcal X\subseteq{\mathbb{R}}^d$ is bounded, and
the stochastic
kernel $\kappa$ is globally Lipschitz continuous as a map from
$\mathcal
X\times\mathcal X$ to $\mathcal P_1(\mathcal X)$, then, for every
$\tau\in
(0,+\infty)$, for sufficiently small $\eps>0$ and sufficiently large
$n\in{\mathbb{N}}$, it holds
\[
{\mathbb{P}}\bigl(\sup\{\w(\mu^n_t,\mu_{t})\dvtx  t\in[0,\tau] \}
\ge
\eps\bigr)\le\exp(-K\eps^3n) ,
\]
where $K$ is a positive constant depending on $\mathcal X$, $\kappa$  and
$\tau$ only.
\end{longlist}
\end{theorem}

Points (a) of Theorem \ref{theomain} will be proved in Section~\ref
{sectweaksol}, while point (b) will be proved in Section \ref
{sectconcentration}. Additional properties of the solution of the
initial value problem associated to (\ref{ODE}) will be studied in
Section~\ref{sectdensitysol}, while Section~\ref{sectexamples}\vadjust{\eject} will
present an analysis of the behavior of the mean-field dynamics for the
model with heterogeneous influential environment, and for the
bounded-confidence opinion dynamics.

\section{Well-posedness of the measure-valued ODE}\label{sectsolutions}
In this section, we shall first prove point (a) of Theorem \ref
{theomain}, that is, that the initial value problem associated to the
ODE (\ref{ODE}) admits a unique solution.
Then, under further technical assumptions, we shall show that, if the
initial measure $\mu_0$ admits a density, so does the solution $\mu
_t$ at any finite time $t$.

\subsection{Weak solutions}\label{sectweaksol}
To start with, we extend the ODE to the space of signed measures
$\mathcal
M^+(\mathcal X)$.
In order to do this, we need to extend the opera\-tor~$F$ and introduce
another operator $G$ in the following way. For $\mu\in\mathcal
M(\mathcal X)$,
put
%
\begin{equation}\label{Gdef}
F(\mu):=F(\mu^+) , \qquad  G(\mu):=\mu
^+(\mathcal
X)\mu,
\end{equation}
where $\mu=\mu^+-\mu^-$ denotes the Hahn--Jordan
decomposition of $\mu\in\mathcal M(\mathcal X)$.
It is not hard to check that both $F$ and $G$ are locally
Lipschitz continuous with respect to the total variation norm,
that is, for every bounded set $\Theta\subseteq\mathcal M(\mathcal X)$, there
exist nonnegative constants $K_F,K_G$ such that
%
\begin{eqnarray}\label{boundedFG}
\|F(\mu_1)-F(\mu_2)\|&\leq&
K_{F}\|\mu_1-\mu_2\| ,\nonumber
\\[-8pt]
\\[-8pt] \|G(\mu_1)-G(\mu_2)\|&\leq&
K_{G}\|\mu_1-\mu_2\|
\nonumber
\end{eqnarray}
for all $\mu_1,\mu_2\in\Theta$.
Moreover,
%
\begin{equation}\label{massconserv} F(\mu)(\mathcal X)=G(\mu
)(\mathcal
X)=\mu(\mathcal X)^2  \qquad  \forall\mu\in\mathcal M^+(\mathcal
X) .
\end{equation}

In the following, we want to study the well-posedness of initial value
problems associated to the measure-valued ODE
%
\begin{equation}
\frac{ d}{ dt}\mu_t=F(\mu_t)-G(\mu_t) ,\label
{generalODE}
\end{equation}
where (\ref{generalODE}) means that, for every $\varphi\in\mathcal
C_0(\mathcal X)$, the real-valued map $t\mapsto\langle\mu_t,\varphi
\rangle$
is differentiable on ${\mathbb{R}}^+$, and satisfies
$\frac{ d}{ dt}\langle\mu_t,\varphi\rangle=\langle F(\mu
_t),\varphi\rangle
-\langle G(\mu_t),\varphi\rangle,$
for every $t>0$. We shall refer to such a $\{\mu_t\dvtx  t\ge0\}$ as a
weak solution of (\ref{generalODE}).
\begin{proposition}\label{existenceuniquenesslemma} Suppose that
$F,G\dvtx \mathcal M(\mathcal X)\to\mathcal M^+(\mathcal X)$ satisfy
properties~(\ref{boundedFG}), and (\ref{massconserv}). Then, for every
$\mu\in\mathcal M^+(\mathcal X)$, there exists a unique solution
$\{\mu_t\dvtx  t\in{\mathbb{R}}^+\}\subseteq\mathcal M^+(\mathcal X)$ to
(\ref{generalODE}) such that $\mu_0=\mu$. Moreover, $\mu_t(\mathcal
X)=\mu(\mathcal X)$ for every $t\geq0$.
\end{proposition}

\begin{pf} For $\tau\in(0,+\infty)$, let $\mathcal C([0,\tau], \mathcal
M(\mathcal
X))$ be the space of continuous curves in $\mathcal M(\mathcal X)$ equipped
with the sup norm
$\|\{\mu_t\}\|_{\tau} := \sup\{\|\mu_t\| \dvtx t\in[0,\tau
] \}$.
Given a curve $\{\mu_s\}\in\mathcal C([0,\tau], \mathcal M(\mathcal
X))$, and a
bounded measurable function $\varphi\in\mathcal C_0(\mathcal X)$, define
%
\begin{eqnarray}\label{definoper}
\langle\Phi(\{\mu_s\})_t,\varphi
\rangle:=\langle\mu
,\varphi\rangle+\int_0^t\langle
F(\mu_s),\varphi\rangle\,ds-\int_0^t\langle G(\mu_s),\varphi
\rangle\,d
s  \nonumber
\\[-8pt]
\\[-8pt] \eqntext{\forall t\in[0,\tau] .}
\end{eqnarray}
Observe that
(\ref{generalODE}) with the initial condition $\mu_0=\mu$ is
equivalent to
%
\begin{equation}
\langle\mu_t,\varphi\rangle=\langle\Phi(\{\mu_s\}
)_t,\varphi
\rangle \qquad  \forall\varphi\in\mathcal
C_0(\mathcal X),\ t\ge0 .\label{equivODEproblem}
\end{equation}
Notice
that, for every $t\in[0,\tau]$, $\Phi(\{\mu_s\})_t$ can be seen as
the difference of two bounded linear positive functionals on $\mathcal
C_0(\mathcal X)$, so that $\Phi(\{\mu_s\})_t\in\mathcal M(\mathcal X)$.
Moreover, the map $t\mapsto\Phi(\{\mu_s\})_t$ is continuous over
$[0,\tau]$, since
%
\begin{eqnarray}\label{Phitcont}
\|\Phi(\{\mu_s\})_{t+\eps}-\Phi(\{\mu_s\})_t\|
&=&\int_{t}^{t+\eps}\|G(\mu_s)\|\,d
s+\int_{t}^{t+\eps}\|F(\mu_s)\|\,ds\nonumber
\\[-8pt]
\\[-8pt]
&\leq&\eps[\|\{G(\mu_s)\}\|_{\tau}+\|\{F(\mu_s)\}\|_{\tau
} ] .
\nonumber
\end{eqnarray}
Therefore, the operator $\Phi$ takes values in $\mathcal C([0,\tau
],\mathcal
M(\mathcal X))$.
Now, let us consider $\Theta:=\{\nu\in\mathcal M(\mathcal X)\dvtx  \|\nu
\|\leq
2\|\mu\|\}$, let $K_F,K_G$ be the Lipschitz constants relative to
$\Theta$ of $F$, and $G$, respectively. For every $\nu\in\Theta$,
(\ref{massconserv}) and (\ref{boundedFG}), imply that
%
\begin{equation}\label{boundedF} \|F(\nu)\|\leq
\|F(\nu)-F(\mu)\|+\|F(\mu)\|\leq4K_F\|\mu\| .
\end{equation}
Similarly,
%
\begin{equation}\label{boundedG} \|G(\nu)\|\leq4K_G\|\mu\| .
\end{equation}
Define now the set $\mathcal S:=\{\{\mu_t\}\in\mathcal C([0,\tau
],\mathcal M(\mathcal
X))\dvtx  \mu_0=\mu, \mu_t\in\Theta,\ \forall t\in[0,\tau]\}$. For
all $\{\mu_t\}\in S$, using (\ref{boundedF}) and (\ref{boundedG}),
and arguing like in (\ref{Phitcont}), we obtain
%
\begin{equation}\label{PhiS}
\|\Phi(\{\mu_t\})\|_{\tau}\leq(1+4\tau K)\|\mu\| ,
\end{equation}
where
$K:=K_F+K_G$. Moreover, if both $\{\mu_t\}$ and $\{\nu_t\}$ belong
to $\mathcal S$, then,
%
\begin{eqnarray}\label{contraction}
 \|\Phi(\{\mu_t\})- \Phi(\{\nu_t\})\|_{\tau}
&=& \sup_{0\leq t\leq
\tau}\int_{0}^{t} \bigl(\|F(\mu_s)-F(\nu_s)\|+\|G(\nu_s)-G(\mu
_s)\| \bigr)\,ds\nonumber
\\[-8pt]
\\[-8pt]
&\leq&\tau K\|\{\mu_t\}-\{\nu_t\}\|_{\tau} .
\nonumber
\end{eqnarray}
We now assume
to have chosen $\tau\in\,]0,\frac1{4K}]$. Then, by (\ref{PhiS}),
$\Phi(\mathcal S)\subseteq\mathcal S$ and, by (\ref{contraction}),
$\Phi$
is a contraction of $\mathcal S$. Hence, by Banach's fixed point
theorem there exists a unique fixed point of $\Phi$ in $\mathcal S$. As
observed, such a fixed point corresponds to a solution $\{\mu_t\}$
of the ODE (\ref{generalODE}) for $t\in[0,\tau]$, with the initial
condition $\mu_0=\mu$. We now show that indeed $\mu_t\in\mathcal
M^+(\mathcal X)$ for $t\in[0,\tau]$. By contradiction, assume that
there exists $B\in\mathcal B(\mathcal X)$ such that $\mu_t(B)<0$ for some
$t\in[0,\tau]$, and let $t^*:=\sup\{s\in[0,t]\dvtx  \mu_s(B)\ge0\}$.
By continuity, $\mu_{t^*}(B)=0$ while $\mu_{s}(B)<0$ for all $s\in\,
]t^*, t]$. This implies that
\[
F(\mu_s)(B)-G(\mu_s)(B)\ge-\mu_s^+(\mathcal X)\mu_s(B)\geq0
 \qquad
\forall s\in\,]t^*,t] .
\]
But then
\[
\mu_{t}(B)=\int_{t^*}^{t}\bigl(F(\mu_s)(B)-G(\mu_s)(B) \bigr)\,d
s\geq0,
\]
which is a contradiction. Hence, $\mu_t\in\mathcal M^+(\mathcal X)$ for
$t\in[0,\tau]$. Notice moreover that, because of property
(\ref{massconserv}), $\mu_t(\mathcal X)=\mu(\mathcal X)$ for all
$t\in
[0,\tau]$. Finally, a standard induction argument allows one to
extend the existence and uniqueness of the solution to the whole
interval $[0,+\infty)$.
\end{pf}

Notice that, when considering an initial condition $\mu_0\in\mathcal
P(\mathcal
X)$, the solution of (\ref{generalODE}) satisfies $\mu_t\in\mathcal
P(\mathcal
X)$ for all $t$,
thus proving point (a) of Theorem \ref{theomain}.

\subsection{Probability density solutions}
\label{sectdensitysol}
We shall now investigate on the existence of density solutions when
the initial condition $\mu_0$ is absolutely continuous with respect to
Lebesgue's measure.

Given the interaction kernel $\kappa( \cdot| \cdot, \cdot
)$, and a nonnegative measure $\mu$ in $\mathcal M^+(\mathcal X)$, we
put
%
\begin{eqnarray}
\label{nu12}   \kappa_1(\mu)(B|y)&:=&\int\kappa(B|x,y)\,d\mu(x)
,\nonumber
\\[-8pt]
\\[-8pt]  \kappa_2(\mu)(B|x)&:=&\int\kappa(B|x,y)\,d\mu(y)
\nonumber
\end{eqnarray}
for
all $B\in\mathcal B(\mathcal X)$, $x,y\in\mathcal X$. The following result
characterizes regularity properties of the solution of the initial
value problem associated to the ODE (\ref{generalODE}).

\begin{proposition}\label{propfinitetimeAC}
Assume that $\mu_0\LL\lambda$, and that
%
\begin{equation}\label{abscont}
\mu\LL\lambda \quad \Longrightarrow \quad
\kappa_1(\mu)( \cdot|y), \kappa_2(\mu)( \cdot|x)\LL\lambda
 \qquad
\forall x\in\mathcal X ,\ \forall y\in\mathcal X .
\end{equation}
Then,
$\mu_t\LL\lambda$, for all $t\in[0,+\infty)$. Moreover, if there
exists $C\in(0,+\infty)$ such that, for all $\mu\LL\lambda$,
%
\begin{equation}
\label{hypLinftycontrol}
\biggl\|\frac{ d\kappa_2(\mu)( \cdot|x)}{ d\lambda}
\biggr\|_{\infty}\le
C \biggl\|\frac{ d\mu}{ d\lambda} \biggr\|_{\infty}
 \qquad
\forall x\in\mathcal X ,
\end{equation}
then, the density
$f_t= d\mu_t/ d\lambda$ satisfies the estimation
%
\begin{equation}\label{Linftycontrol}
\|f_t\|_{\infty}\le
\|f_0\|_{\infty}e^{Ct}
 \qquad
\forall t\in[0,+\infty) .
\end{equation}
\end{proposition}

\begin{pf} For every finite time $t\in[0,+\infty)$, consider Lebesgue's
decomposition $\mu_t=\mu_t^{a}+\mu_t^{s}$, where
$\mu_t^{a}\LL\lambda$, and $\mu_t^{s}$ and $\lambda$ are singular.
It follows from (\ref{abscont}) that,
$\kappa_2(\mu_t^a)( \cdot|x)\LL\lambda$ for all $x\in\mathcal X$.
Then, for any $B\in\mathcal B(\mathcal X)$ such that $\lambda(B)=0$,
one has
\[
\int\!\!\int\kappa(B|x,y)\,d\mu_t^a(x)\,d\mu_t(y)=\int\,d
\kappa_2(\mu_t^a)(B|x)\,d\mu_t(x)=0 .
\]
Similarly, one can show that $\int\!\!\int\kappa(B|x,y)\,d\mu
^s_t(x)\,d\mu^a_t(y)=0$.
Hence,
\begin{eqnarray*}
F(\mu_t)(B)&=&\int\!\!\int\kappa(B|x,y)\,d\mu
_t(x)\,d\mu_t(y)
\\&=&\int\!\!\int\kappa(B|x,y)\,d\mu_t^a(x)\,d\mu_t(y)+
\int\!\!\int\kappa(B|x,y)\,d\mu_t^s(x)\,d\mu_t^a(y)\\
&&{}+\int\!\!
\int\kappa(B|x,y)\,d\mu_t^s(x)\,d\mu_t^s(y)
\\&=&\int\!\!\int\kappa(B|x,y)\,d\mu_t^s(x)\,d\mu_t^s(y)
\\&=&F(\mu_t^s)(B)
\end{eqnarray*}
for all $B\in\mathcal B(\mathcal X)$ such that $\lambda(B)=0$.
This readily implies that $\mu_t^{s}$ satisfies
\[
\frac{ d}{ dt}\mu^{s}_t=F(\mu_t^{s})-\mu^s_t .
\]
Since $\mu_0^{s}=0$ by assumption, it follows that $\mu_t^{s}=0$ for
all $t\ge0$.

Assume now that (\ref{hypLinftycontrol}) holds true. For any
$\varphi\in\mathcal C_c(\mathcal X)$, H\"older's inequality, and
(\ref{hypLinftycontrol}) imply that
\begin{eqnarray*}
\langle F(\mu_t),\varphi\rangle
&=&
\int\!\!\int\!\!\int\varphi(z)\,d\kappa(z|x,y)\,d\mu
_t(x)\,d\mu_t(y)\\
&=&
\int\!\!\int\varphi(z)\frac{ d\kappa_2(\mu_t)(z|x)}{ d
\lambda}\,d\lambda(z) \,d\mu_t(x)\\
&\le&
\int\biggl\|\frac{ d\kappa_2(\mu_t)(z|x)}{ d\lambda
} \biggr\|_{\infty}\|\varphi\|_1\,d\mu_t(x)\\
&\le&C\|f_t\|_{\infty}\|\varphi\|_1 .
\end{eqnarray*}
It follows that, for
all nonnegative-valued $\varphi\in\mathcal C_c(\mathcal X)$,
\begin{eqnarray*} \int\varphi(x)f_t(x)\,d\lambda(x)&=&\int
\varphi(x)f_0(x)\,d\lambda(x)+
\int_0^t\bigl(\langle F(\mu_s),\varphi\rangle-\langle\mu_s,\varphi
\rangle
\bigr) \,ds\\
&\le&\|f_0\|_{\infty}\|\varphi\|_1+\int_0^t\langle F(\mu
_s),\varphi\rangle\, ds\\
&\le&\|\varphi\|_1\biggl(\|f_0\|_{\infty}+C\int_0^t\|f_s\|_{\infty
} \,d
s \biggr) .
\end{eqnarray*}
Then, by the isometry of $L^{\infty}(\mathcal X)$ with the dual of
$L^1(\mathcal X)$, the fact that $f_t$ is nonnegative valued, and the
density of $\mathcal C_c(\mathcal X)$ in $L^1(\mathcal X)$, one gets that
\begin{eqnarray*}
\|f_t\|_{\infty}&=&
\sup\biggl\{\int\varphi(x)f_t(x)\,dx\dvtx  \varphi\in L^1(\mathcal X),\
\|\varphi\|_1\le1 \biggr\}\\
&=&\sup\biggl\{\int\varphi(x)f_t(x)\,dx\dvtx  \varphi\in\mathcal
C_c(\mathcal
X), \varphi\ge0, \|\varphi\|_1\le1 \biggr\}\\
&\le&\|f_0(x)\|_{\infty}+C\int_0^t\|f_s\|_{\infty} \,ds .
\end{eqnarray*}
By Gronwall's lemma, this readily implies
(\ref{Linftycontrol}).\vspace*{-2pt}
\end{pf}

The technical condition on the stochastic kernel $\kappa$ is
actually verified in many important cases encompassing the bounded
confidence dynamics (\ref{boundedconfidence}) as well as the
Gaussian interaction model (\ref{Gaussiandecay}).\vspace*{-2pt}

\begin{corollary}
Assume that the interaction kernel $\kappa$ is the form (\ref
{interacting kernels}) with
$\theta^i( \cdot|x,y)=\delta_{\omega(|x-y|)}$ and
$\theta^e( \cdot|x,z)=\delta_{\upsilon(|x-z|)}$ where
$\omega\dvtx {\mathbb{R}}^+\to[0,\omega_0]$, $\omega_0\in[0,1[$, and
$\upsilon\dvtx {\mathbb{R}}^+\to[0,\upsilon_0]$, $\upsilon_0\in[0,1[$,
are both
nonincreasing and piecewise~$\mathcal C^1$. If $\mu_0\LL\lambda$, then
$\mu_t\LL\lambda$, for all $t\in[0,+\infty)$ and the relative
densities satisfy condition (\ref{Linftycontrol}).\vspace*{-2pt}
\end{corollary}

\begin{pf} We shall show that the conditions of Proposition
\ref{propfinitetimeAC} are satisfied in this case. Fix $y\in\mathcal X$
and consider the function
$x\mapsto\ov\omega(|x-y|)x+\omega(|x-y|)y$. The assumption on
$\omega$ ensures that it is an invertible transformation in~$x$
and a simple geometric consideration shows that the inverse has
the form\looseness=-1
\[
x=g(w,y)=y+\alpha(|w-y|)(w-y),
\]
where $\alpha\dvtx {\mathbb{R}}^+\to
{\mathbb{R}}^+$ is such that $\alpha(\ov\omega(t)t)\ov\omega
(t)=1$ for all
$t\geq0$. The function $\eta(t)=\ov\omega(t)t$ is strictly
increasing, hence invertible and we can thus write
$\alpha(s)=[\ov\omega(\eta^{-1}(s))]^{-1}$. $\alpha$ is thus also a
piecewise $C^1$ function as well as $g(\cdot, y)$ whose Jacobian can
easily be shown to be
\[
D_wg(w,y)=\alpha(|w-y|)I+\nabla\alpha(|w-y|)(w-y)^t.
\]
Straightforward computation show that $D_wg$ is bounded in the
pair $(w,y)$. Similarly, the function
$x\mapsto\ov\upsilon(|x-z|)x+\omega(|x-z|)z$ admits an inverse in $x$,
$x=h(w,z)$ whose Jacobian $D_wh$ is bounded in the pair $(w,z)$.


Then, if $\mu$ is absolutely continuous with
density $f$, one has for all nonnegative real-valued $\varphi\in
L^1(\mathcal X)$, and $y\in\mathcal X$,
\begin{eqnarray*}\langle\kappa_1(\mu)( \cdot|y),\varphi\rangle
&=&\alpha\int\varphi\bigl(\ov\omega(|x-y|)x+\omega
(|x-y|)y \bigr)f(x)\,d\lambda(x)
\\&&{}+\ov\alpha\int\!\!\int\varphi\bigl(\ov\upsilon(|x-z|)x+\upsilon
(|x-z|)z\bigr)f(x)\,d\lambda(x)\,d\psi(z)\\
&=&\alpha\int\varphi(w )f(g(w,y))|D_wg(w,y)|\,d
\lambda(w)\\
&&{}+\ov\alpha\int\!\!\int\varphi(w
)f(h(w,z))|D_wh(w,z)|\,d\lambda(w)\,d\psi(z)\\
&\le&C_1\|f\|_{\infty}\|\varphi\|_1 ,
\end{eqnarray*}
where
$C_1:=\alpha\|D_wg(w ,y)\|_{\infty}+\ov\alpha\|D_wh(w
,z)\|_{\infty}$.
Similarly, one shows that there exists some constant $C_2>0$ such
that $\langle\kappa_2(\mu)( \cdot|x),\varphi\rangle\le
C_2\|f\|_{\infty}\|\varphi\|_1$, for all nonnegative valued
$\varphi\in L_1(\mathcal X)$  and $x\in\mathcal X$. As a consequence,
$\nu_i( \cdot|y)\LL\lambda$, for all $y\in\mathcal X$  and $i=1,2$,
and (\ref{hypLinftycontrol}) holds. Therefore, the claim follows from
Proposition~\ref{propfinitetimeAC}.
\end{pf}


\section{Behavior of the mean-field dynamics}\label{sectexamples}
This section is devoted to a~dee\-per analysis of the ODE
(\ref{generalODE}) for the state-independent gossip model with
heterogeneous influential environment, and the bounded-confidence
opinion dynamics, respectively. In particular, we shall
investigate the limit behavior as $t$ grows large, showing that,
in both models, $\mu_t$ converges weakly to an asymptotic opinion
measure. The behavior of the two models, and their analysis,
however, differ substantially. For the state-independent gossip
model with heterogeneous influential \mbox{environment}, the ODE governing
the mean-field dynamics is linear, and can be analyzed by
iteratively solving the lower-triangular linear system of ODEs
governing the various moments behavior. In this case, the
asymptotic opinion measure is independent of the initial value, it
is characterized by its moments, and is absolutely continuous if
so is the influential environment. In fact, one could show that
the corresponding finite population Markov process is ergodic. In
contrast, the ODE governing the mean-field dynamics of the bounded
confidence model is nonlinear, and convergence is shown by a
Lyapunov argument. The asymptotic opinion measure is given by a
convex combination of deltas, whose number and position typically
depends on the initial condition. Indeed, the corresponding finite
population Markov process is typically not ergodic, in this
case.\looseness=-1

\subsection{Gossip model with heterogeneous influential
environment}\label{sectinfluential}
We start by analyzing the case when the stochastic kernel
$\kappa( \cdot| \cdot, \cdot)$ has the form
(\ref{interacting kernels}), with constant weights:
$\theta^i( \cdot|x,y)=\delta_{\omega}( \cdot)$,
$\theta^e( \cdot|x,y)=\delta_{\upsilon}( \cdot)$, for some fixed
$\omega,\upsilon\in[0,1]$. Throughout this section, we shall
assume an exponential bound on the moments of both $\mu_0$ and
$\psi$, that is,
%
\begin{equation}\label{expboundmoments}
\sup_{k\in{\mathbb{N}}}\biggl(\int|x|^k\,d\mu_0(x)
\biggr)^{1/k}<+\infty
, \qquad
\sup_{k\in{\mathbb{N}}}\biggl(\int|x|^k\,d\psi(x) \biggr)^{1/k}<+\infty
.
\end{equation}
Clearly, (\ref{expboundmoments}) is automatically satisfied when
$\mathcal X$ is bounded. Let us fix some $z\in{\mathbb{R}}^d$, and
consider the
$z$-weighted moments of $\mu_t$ and $\psi$, respectively,
\[
m^{(k)}_t:=\int(x\cdot z)^k\,d\mu_t(x) , \qquad  n^{(k)}_t:=\int
(x\cdot y)^k\,d\psi(y) , \qquad  k\in{\mathbb{Z}}^+ .
\]
%
The following result characterizes their evolution in time.
%
\begin{proposition} The $z$-weighted moments satisfy
%
\begin{eqnarray}
\frac{ d}{ d
t}m^{(1)}_t&=&{\ov\alpha}\upsilon\bigl(n^{(1)}-m_t^{(1)} \bigr) ,\label
{ODE1}
\\
\label{ODEk}\frac{ d}{ d
t}m^{(k)}_t&=&-\gamma_km_t^{(k)}+f_k\bigl(m_t^{(1)},\ldots
,m_t^{(k-1)} \bigr)+{\ov\alpha}\upsilon^kn^{(k)} , \qquad  k\ge2 ,
\end{eqnarray}
where
\begin{eqnarray*}
 \gamma_k&:=&1-\alpha(\ov{\omega} ^k+\omega^k )-{\ov\alpha
}{\ov\upsilon}^k ,
\\
  f_k\bigl(m_t^{(1)} ,\ldots,m_t^{(k-1)} \bigr)&:=& \sum_{j=1}^{k-1}
\binom
kj \bigl( \alpha\ov{\omega}^j\omega^{k-j}m_t^{(j)}m_t^{(k-j)}+{\ov
\alpha}{\ov\upsilon}^j\upsilon^{k-j}m_t^{(j)}n^{(k-j)} \bigr) .
\end{eqnarray*}
\end{proposition}

\begin{pf}
For the first moment, one has
\begin{eqnarray*}\frac{ d}{ dt}m^{(1)}_t
& = & \alpha\int\!\!\int\bigl((\ov\omega
x+\omega y)\cdot z \bigr)\,d\mu_t(x)\,d\mu_t(y)\\
&&{}+\ov\alpha\int\!\!\int\bigl((\ov\upsilon x+\upsilon y)\cdot
z \bigr)\,d\mu_t(x)\,d\psi(y)-m^{(1)}_t\\
& = &\ov\alpha\upsilon n^{(1)}-\ov\alpha\upsilon
m^{(1)}_t ,
\end{eqnarray*}
which proves (\ref{ODE1}). For $k\ge2$, one has
\begin{eqnarray*}
&&\int(\ov\omega x\cdot z+\omega y\cdot z )^k\,d
\mu_t(x)\,d\mu_t(y)\\
&& \qquad =(\ov\omega^k+\omega^k)\int(x\cdot z)^k\,d\mu
_t(x)\\
&& \qquad \quad  {}+\sum_{j=1}^{k-1}\binom{k}{j}\ov{\omega}^j\omega^{k-j}\int
(x\cdot z)^j\,d\mu_t(x)\int(y\cdot z)^{k-j}\,d\mu_t(y)\\
&& \qquad =(\ov\omega^k+\omega^k)m^{(k)}_t+
\sum_{j=1}^{k-1}\binom kj\ov{\omega}^j\omega
^{k-j}m_t^{(j)}m_t^{(k-j)} ,
\end{eqnarray*}
and, similarly,
\begin{eqnarray*}
&&\int\!\!\int \bigl((\ov\upsilon x+\upsilon y) \cdot z
\bigr)^k\,d
\mu_t(x)\,d\psi(y)\\
&& \qquad  =
\ov\upsilon^km^{(k)}_t+\sum_{j=1}^{k-1}\binom kj {\ov\upsilon
}^j\upsilon
^{k-j}m_t^{(j)}n^{(k-j)}+\upsilon^kn^{(k)}.
\end{eqnarray*}
From the two identities above, it follows that
\begin{eqnarray*}
\frac{ d}{ dt}m^{(k)}_t
&=&\alpha\int\!\!\int\bigl((\ov\omega x+\omega y)\cdot z
\bigr)^k\,d\mu_t(x)\,d\mu_t(y)\\
&&{}+\ov\alpha\int\!\!\int\bigl((\ov\upsilon x+\upsilon y)\cdot z
\bigr)^k\,d\mu_t(x)\,d\psi(y)-m^{(k)}_t\\
&=&-\gamma_km_t^{(k)}+f_k\bigl(m_t^{(1)},\ldots,m_t^{(k-1)}
\bigr)+{\ov\alpha}\upsilon^kn^{(k)} ,
\end{eqnarray*}
which proves (\ref{ODEk}).
\end{pf}

\begin{example}
In the special case when $\alpha=1$, namely when there is no
influential environment, we obtain from (\ref{ODE1}) that
$\frac{ d}{ dt}\int x\,d\mu_t(x)=0$, so that the first moment is
constant. On the other hand, the variance
\[
v_t:=\int\biggl|x-\int y\,d\mu_0(y) \biggr|^2\,d\mu_t(x)
\]
satisfies
$\frac{ d}{ dt}v_t=-2\omega\ov{\omega} v_t$.
Hence,
\[
v_t=v_0e^{-\omega\ov{\omega} t} ,
\]
that is, $\mu_t$ converges to a delta centered in the average initial
opinion exponentially fast in $t$.
\end{example}

We now focus on the limit as $t\to+\infty$ for the general case.
An inductive argument proves the following result.
\begin{lemma}
Assume $\alpha<1$. Then,
for every $z\in{\mathbb{R}}^d$, the $z$-weighted moments of $\mu_t$ satisfy
%
\begin{equation}\label{inductionmoments}
\lim_{t\to\infty}m^{(k)}_t=m^{(k)}_{\infty} , \qquad  k\in{\mathbb
{Z}}^+
,
\end{equation}
where $m^{(k)}_{\infty}$ can be recursively evaluated by
%
\begin{equation}\label{momentcomp}m^{(1)}_{\infty}:=n^{(1)} , \qquad
m^{(k+1)}_{\infty}=\gamma_{k+1}^{-1}\bigl[f_{k+1}\bigl(m^{(1)}_{\infty
},\ldots,m^{(k)}_{\infty} \bigr)+{\ov\alpha}\upsilon
^kn^{(k+1)}
\bigr] .
\end{equation}
\end{lemma}

\begin{pf} For $k=1$, the solution of the ODE (\ref{ODE1}) is easily
found to be
%
\begin{equation}\label{ODE1sol}m_t^{(1)}=e^{-{\ov\alpha}\upsilon t}m^{(1)}_0+
(1-e^{-{\ov\alpha}\upsilon t} )n^{(1)} ,
\end{equation}
so that equation (\ref{inductionmoments}) clearly holds. Moreover,
assume that equation (\ref{inductionmoments}) holds for every $k\in\{
1,\ldots,j-1\}$, and define
$\chi_t^{(j)}:=f_{j}(m_t^{(1)},\ldots,m_t^{(j-1)} )$ for
$t\in[0,+\infty]$. Then, the continuity of $f_j$ implies that
$\lim_{t\to\infty}\chi_t^{(j)}=\chi_{\infty}^{(j)}$.
Solving the ODE (\ref{ODEk}) gives
%
\begin{equation}\label{ODEksol}
m^{(j)}_t=\int_0^te^{-\gamma_{j}(t-s)}\bigl(\chi_t^{(j)}+{\ov\alpha
}\upsilon^jn^{(j)} \bigr)\,d
s+e^{-\gamma_jt}m_0^{(j)} .
\end{equation}
Clearly, the second addend of the
right-hand side of (\ref{ODEksol}) converges to zero for
$t\to\infty$. On the other hand, the convergence of $\chi_t^{(j)}$
implies that
\begin{eqnarray*}
\lim_{t\to\infty}\int_0^te^{-\gamma
_{j}(t-s)}\bigl(\chi_t^{(j)}+{\ov\alpha}\upsilon^jn^{(j)} \bigr)\,ds&=&
\bigl(\chi_{\infty}^{(j)}+{\ov\alpha}\upsilon^jn^{(j)}
\bigr)\lim_{t\to\infty}\int_0^t e^{-\gamma_{j}(t-s)}\,ds\\
&=&\gamma_j^{-1}\bigl(\chi_{\infty}^{(j)}+{\ov\alpha}\upsilon
^jn^{(j)} \bigr) .
\end{eqnarray*}
The foregoing, together with (\ref{ODEksol}), implies the claim.
\end{pf}

We are now in a position to prove the following result for the
convergence of $\mu_t$.
\begin{proposition}\label{propinfluentialconvergence}
Assume that (\ref{expboundmoments}) holds.
Then
\[
\lim_{t\to\infty}\mu_t=\mu_{\infty} ,
\]
weakly, where $\mu_{\infty}\in\mathcal P(\mathcal X)$ is uniquely
characterized by its moments $m^{(k)}_{\infty}$.
\end{proposition}

\begin{pf}
It follows from (\ref{expboundmoments}) that there exists some finite
$M\in{\mathbb{R}}^+$ such that
%
\begin{equation}\label{m<=M^k}
\bigl|m_0^{(k)}\bigr|\le|z|^kM^k , \qquad
\bigl|n^{(k)}\bigr|\le
|z|^kM^k
\end{equation}
for all $z\in{\mathbb{R}}^d$  and $k\in{\mathbb{N}}$.
Now, an inductive argument shows that
%
\begin{equation}\label{mt<=M^k}
\bigl|m_t^{(k)}\bigr|\le|z|^kM^k  \qquad
\forall t\in
[0,+\infty] ,\ z\in{\mathbb{R}}^d
\end{equation}
for all $k\in{\mathbb{N}}$. In fact, (\ref{ODE1sol}) and (\ref{m<=M^k})
immediately imply that (\ref{mt<=M^k}) holds for $k=1$.
Moreover, if (\ref{mt<=M^k}) holds for all $k\in\{1,\ldots,j-1\}$,
then (\ref{ODEksol})  and (\ref{mt<=M^k}) give
\begin{eqnarray*}
\bigl|m^{(j)}_t\bigr| &\le&
\int_0^te^{-\gamma_{j}(t-s)}\bigl(\bigl|f_{j}\bigl(m_t^{(1)},\ldots
,m_t^{(j-1)} \bigr) \bigr|+{\ov\alpha}\upsilon^j\bigl|n^{(j)}\bigr|
\bigr)\,d
s+e^{-\gamma_jt}\bigl|m_0^{(j)}\bigr|\\
&\le& \int_0^te^{-\gamma_{j}(t-s)}M^j|z|^j\gamma_j\,d
s+e^{-\gamma_jt}M^j|z|^j\\
&=& M^j|z|^j .
\end{eqnarray*}

Let us consider the characteristic functions
$\phi_t(z):=\int\exp(iz\cdot x)\,d\mu_t(x)$  and, for $k\in
{\mathbb{Z}}
^+$, define $a_t(k):=i^km_t^{(k)}/k!$, $b(k):=M^k|z|^k/k!$, and observe
that $\sum_{k\in{\mathbb{Z}}^+}b(k)=\exp(M|z|)$. One has that
\[
\phi_t(z)=\int\sum_{k\in{\mathbb{Z}}^+}\frac{(iz\cdot
x)^k}{k!}\,d\mu_t(x)=
\sum_{k\in{\mathbb{Z}}^+}\frac{i^k}{k!}\int(x\cdot z)^k\,d\mu_t(x)=
\sum_{k\in{\mathbb{Z}}^+}a_t(k) ,
\]
where the exchange between the series and the integral is justified by
Lebes\-gue's dominated convergence theorem, since
\[
\biggl|\sum_{0\le k\le n}\frac{i^k}{k!}(x\cdot z)^k \biggr|\le
\sum_{0\le k\le n}b(k)\le\exp(M|z|) .
\]
Moreover, observe that, since $|a_t(k)|\le b(k)$, another application
of Lebes\-gue's dominated convergence theorem gives
\[
\lim_{t\to\infty}\phi_t(z)=
\lim_{t\to\infty}\sum_{k\in{\mathbb{Z}}^+}a_t(k)=
\sum_{k\in{\mathbb{Z}}^+}a_{\infty}(k)=:\phi_{\infty}(z) .
\]
Hence, $\phi_{t}(z)$ converges pointwise to $\phi_{\infty}(z)$,
which in turn can be easily verified to be continuous at $0$.
Then, the claim follows from L\'evy's continuity theorem (\cite{Borkar}, Theorem~2.5.1).
\end{pf}

Observe that, for all $\alpha\in(0,1)$, the limit measure $\mu
_{\infty}$ is independent of the initial condition $\mu_0$, and
depends only on the influential environment~$\psi$, as well as on the
parameters $\alpha$, $\omega$  and $\upsilon$.
Notice that the first moment satisfies $m^{(1)}_{\infty}=n^{(1)}$.
In contrast, if $\psi\ne\delta_{x_0}$, it easily seen that
$m^{(k)}_{\infty}\ne n^{(k)}$ for $k\ge2$, so that in particular $\mu
_{\infty}\ne\psi$.
On the other hand, it follows from (\ref{momentcomp}) that, if $\psi
\ne\delta_{x_0}$, then the variance of $\mu_{\infty}$ is positive,
so that $\mu_{\infty}\ne\delta_{x_0}$. This result may be
interpreted as showing that the presence of an heterogeneous
influential environment prevents the population from achieving an
asymptotic opinion agreement. In fact, as shown in the following
proposition, the asymptotic opinion distribution $\mu_{\infty}$ is
absolutely continuous whenever so is the influential environment $\psi$.
\begin{proposition}\label{propinfluentiallimitAC}
Assume $\psi\LL\lambda$. Then $\mu_{\infty}\LL\lambda$ for all
$\alpha\in[0,1)$.
\end{proposition}

\begin{pf}
For $\mu,\nu\in\mathcal P(\mathcal X)$, $\gamma\in[0,1]$, define
$\ov\gamma
:=1-\gamma$, and
\[
L_{\gamma}(\mu,\nu)\in\mathcal P(\mathcal X), \qquad  \langle
L_{\gamma}(\mu,\nu
),\varphi\rangle=
\int\!\!\int\varphi(\ov\gamma x+\gamma y)\,d\mu(x)\,d\nu(y)
\]
for every $\varphi\in\mathcal C_b(\mathcal X)$.
Since $L_{\gamma}$ is a rescaled convolution operator, and since $\psi
\LL\lambda$, one has that $L_{\upsilon}(\mu,\psi)\LL\lambda$.
Similarly, $L_{\omega}(\mu_{\infty},\mu_{\infty})=\alpha L(\mu
^s_{\infty},\mu^s_{\infty})$, where $\mu^s_{\infty}$ is the
singular part of $\mu_{\infty}$.
Combining this with the fact that the asymptotic measure satisfies
\[
\mu_{\infty}=F(\mu_{\infty})=\alpha L_{\omega}(\mu,\mu)+\ov
\alpha L_{\upsilon}(\mu,\psi) ,
\]
one gets that
\[
\mu^s_{\infty}(\mathcal X)=\alpha(L_{\omega}(\mu^s_{\infty},\mu
^s_{\infty}) )(\mathcal X)=\alpha(\mu^s(\mathcal X) )^2 .
\]
Therefore,
\[
\mu^s_{\infty}(\mathcal X)\bigl(1-\alpha\mu^s_{\infty}(\mathcal
X) \bigr)=0 .
\]
Since $\mu^s_{\infty}(\mathcal X)\le1$ and $\alpha<1$, this
necessarily implies that $\mu^s_{\infty}(\mathcal X)=0$.
\end{pf}

\begin{figure}

\includegraphics{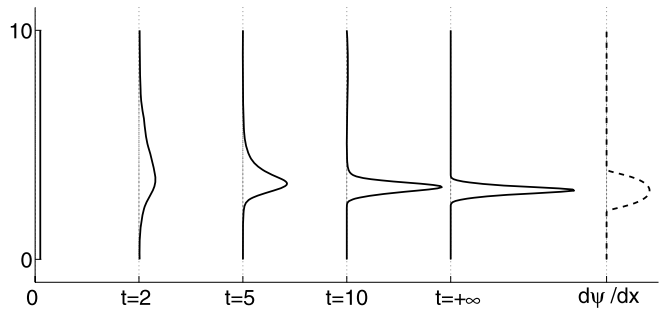}

\caption{Behavior in time of the ODE solution
in $d=1$, with
initial condition $\mu_0$ uniform over $(0,10)$, heterogeneous
environment $ d\psi(x)=\exp(-(1-(x-3)^2)^{-1})\mathbh
{1}_{(2,4)}(x)\,d
x$, and parameters $\alpha=0.5$, $\omega=0.5$  and $\upsilon=0.5$. The
Radon--Nikodym derivates of the asymptotic measure $\mu_{\infty}$,
and of the influential environment $\psi$ (dashed) are plotted as
a reference.}\label{figinfluential}
\end{figure}

Figure~\ref{figinfluential} reports numerical simulations of the
mean-field dynamics, when started from a uniform distribution over an
interval, and influenced by an absolutely continuous environment.
Coherently with Proposition \ref{propfinitetimeAC}, the solution
remains absolutely continuous during its evolution. As $t$ grows large,
$\mu_t$ converges to a limit measure whose first moment coincides with
that of~$\psi$, and which is absolutely continuous, as predicted by
Propositions \ref{propinfluentialconvergence}, and~\ref
{propinfluentiallimitAC}, respectively. Such a limit density may be
interpreted as resulting from a~tension between the aggregating forces
represented by the first addend in the right-hand side of (\ref
{interacting kernels}), and the environment's influence captured by the
second addend in the right-hand side of~(\ref{interacting kernels}).

\subsection{Bounded confidence opinion dynamics}\label{sectboundedconfidence}
We analyze now the case when $\kappa( \cdot| \cdot, \cdot
)$ is in the form (\ref{interacting kernels}) with $\alpha=1$, and
weight distribution $\theta( \cdot|x,y):=\theta^i( \cdot
|x,y)$ supported on $[0,\omega_0]$ for some $\omega_0\in[0,1[$, and
satisfying the symmetry assumption
%
\begin{equation}\label{symmetry}
\theta( \cdot
|x,y)=\theta( \cdot|y,x)
\end{equation}
for all $x,y\in\mathcal X$. The
following result states weak convergence of $\mu_t$.

\begin{proposition}\label{propboundedconfidenceconvergence}
Assume that $\int|x|^2\,d\mu_0(x)<\infty$. Then, there exists $\mu
_{\infty}\in\mathcal P(\mathcal X)$ such that
\[
\lim_{t\to\infty}\mu_t=\mu_{\infty} ,
\]
weakly.
\end{proposition}

\begin{pf} We start by proving that the second moment
$m^{(2)}_t\!:=\!\int\!|x|^2\,d\mu_t(x)$ is a Lyapunov function for the system. Observe
that, for all
$x,y\in{\mathbb{R}}^d$, $\omega\in[0,1]$, $\bar\omega=1-\omega$,
one has
\[
|x+\omega(y-x) |^2+|y+\omega(x-y) |^2=
(\ov{\omega}
^2+\omega^2 )(|x|^2+|y^2| )+4\omega\ov{\omega}
x\cdot y ,
\]
so that
\begin{eqnarray*}
2\omega\ov{\omega} |x-y|^2&=&
2\omega\ov{\omega} (|x|^2+|y|^2-2x\cdot y )\\
&=&(1-\omega^2-\ov{\omega} ^2 )(|x|^2+|y|^2
)-4\omega\ov{\omega} x\cdot y\\
&=&|x|^2+|y^2|-|x+\omega(y-x) |^2-|y+\omega(x-y)
|^2 .
\end{eqnarray*}
From the foregoing, and the symmetry of $\theta( \cdot|x,y)$, it
follows that
%
\begin{eqnarray}\label{deVtdet}
\frac{ d}{ dt}m^{(2)}_t& = &\int|x|^2\,dF(\mu
_t)(x)-m^{(2)}_t\nonumber\\
& = &\int\!\!\int\!\!\int\bigl(|x+\omega
(y-x) |^2-|x|^2 \bigr)\,d\theta(\omega|x,y)\,d\mu_t(x)\,d
\mu_t(y)\nonumber\\
& =& \frac12\int\!\!\int\!\!\int \bigl(
|x+\omega(y-x) |^2
\nonumber
\\[-8pt]
\\[-8pt]&&\hphantom{\frac12\int\!\!\int\!\!\int \,}
{}+ |y+\omega(x-y) |^2 -|x|^2-|y|^2
\bigr)\,d\theta(\omega|x,y)\,d\mu_t(x)\,d\mu_t(y)\nonumber\\
& = &-\int\!\!\int\!\!\int
\omega(1-\omega)|x-y|^2\,d\theta(\omega|x,y)\,d\mu
_t(x)\,d\mu_t(y) \nonumber\\
& \leq &-(1-\omega_0)\Upsilon_t ,
\nonumber
\end{eqnarray}
where
\[
\Upsilon_t:= \int\!\!\int\!\!\int
\omega|x-y|^2\,d\theta(\omega|x,y)\,d\mu_t(x)\,d\mu_t(y) .
\]
Hence, in particular, $\frac{ d}{ dt}m^{(2)}_t\le0$, so that $m^{(2)}_t$
is nonincreasing, and therefore convergent. Define
$m^{(2)}_{\infty}:=\lim_{t\to\infty}m^{(2)}_t$ and observe
that (\ref{deVtdet}) implies that
%
\begin{eqnarray}\label{Cauchy}
\lim_{t\to\infty}\int_0^t\Upsilon_s\,ds
&\leq&\lim_{t\to\infty}-\frac1{1-\omega_0 }\int_0^t\frac
{ d}{ ds}m^{(2)}_s\,ds\nonumber\\
&=&\lim_{t\to\infty}\frac{m^{(2)}_0-m^{(2)}_t}{1-\omega_0 }\\
&=&\frac{m^{(2)}_0-m^{(2)}_{\infty}}{1-\omega_0 } .\nonumber
\end{eqnarray}
Now, for any smooth and compact-supported test function
$\varphi\in\mathcal\mathcal C^{\infty}_{c}({\mathbb{R}}^d)$, we
can write
%
\begin{equation}\label{Taylor}
\varphi\bigl(x+\omega(y-x)\bigr)-\varphi
(x)=\omega
(y-x)\cdot
\nabla\varphi(x)+r(x,y) ,
\end{equation}
with $|r(x,y)|\leq
\omega^2|y-x|^2\Phi$ where $\Phi:=\|D^2\varphi\|$. Moreover, again
from the symmetry of $\theta( \cdot|x,y)$, one has
%
\begin{eqnarray}\label{gradestim}
 &&\biggl|\int\!\!\int\!\!\int
\omega(y-x)\cdot\nabla\varphi(x)\,d\theta(\omega|x,y)\,d\mu
_t(x)\,d\mu_t(y) \biggr|\nonumber\\
&& \qquad =
\frac12\biggl|\int\!\!\int\!\!\int\omega(y-x)\cdot\bigl(\nabla\varphi
(x)-\nabla\varphi(y) \bigr)\,d\theta(\omega|x,y)\,d\mu_t(x)\,d
\mu_t(y) \biggr|\nonumber\\
&& \qquad \le
\frac12\int\!\!\int\!\!\int\omega|y-x|\bigl|\nabla\varphi
(x)-\nabla\varphi(y) \bigr|\,d\theta(\omega|x,y)\,d\mu_t(x)\,d
\mu_t(y)\\
&& \qquad \le
\frac{\Phi}2\int\!\!\int\!\!\int\omega|x-y|^2\,d\theta(\omega
|x,y)\,d\mu_t(x)\,d\mu_t(y)\nonumber\\
&& \qquad \le
\frac{\Phi}2\Upsilon_t\nonumber
.
\end{eqnarray}

From (\ref{Taylor}) and (\ref{gradestim}), it follows that
\begin{eqnarray*}
|\langle F(\mu_t)-\mu_t,\varphi\rangle
|
&=&
\biggl|\int\!\!\int\!\!\int \bigl(\varphi\bigl(x+\omega(y-x)
\bigr)-\varphi(x) \bigr)\,d\theta(\omega|x,y)\,d\mu_t(x)\,d\mu
_t(y) \biggr|\\
&\le& \biggl|\int\!\!\int\!\!\int\omega(y-x)\cdot\nabla
\varphi(x)\,d\theta(\omega|x,y)\,d\mu_t(x)\,d\mu_t(y) \biggr|\\
&&{} +\Phi\int\!\!\int\!\!\int\omega^2|x-y|^2\,d\theta
(\omega|x,y)\,d\mu_t(x)\,d\mu_t(y)\\
&\le&
\frac{3\Phi}{2}\Upsilon_t ,
\end{eqnarray*}
so that
\begin{eqnarray*}
\lim_{t\to\infty}\int_0^t|\langle F(\mu_s)-\mu_s,\varphi\rangle
|\,ds
&\le&
\frac{3\Phi}{2}\lim_{t\to\infty}\int_0^t\Upsilon_s\,ds\\
&\le& \frac{3\Phi}{2(1-\omega_0)}\bigl(m^{(2)}_0-m^{(2)}_{\infty}\bigr)
.
\end{eqnarray*}
Therefore, in particular, the limit
\[
\lim_{t\to\infty}\langle\mu_t,\varphi\rangle= \lim_{t\to\infty
}\int
_0^t\langle F(\mu_s)-\mu_s,\varphi\rangle\,ds
\]
exists and is finite. From the arbitrariness of $\varphi\in\mathcal
C^{\infty}_c(\mathcal X)$, it follows that $\mu_t$ converges in the sense
of distributions.
Finally, notice that, since the second moment is bounded, the family $\{
\mu_t\dvtx  t\in{\mathbb{R}}^+\}$ is tight, hence $\mu_t$ converges in
$\mathcal P(\mathcal
X)$.
\end{pf}

If we make the further assumption that the weight
$\omega\sim\theta( \cdot|x,y)$ is almost surely strictly
positive in a
neighborhood of the diagonal $\{(x,x)\dvtx  x\in\mathcal X\}$, we have the
following characterization of the equilibrium points.
\begin{proposition}\label{propboundedconfidenceasymptotic}
Let $R>0$ be such that
\[
\delta(R):=\inf\{\omega\dvtx  \operatorname{supp}(\theta( \cdot
|x,y))\subseteq
[\omega,1]\ \forall x,y\in\mathcal X, |x-y|<R\}>0.
\]
Then $\mu_{\infty}$ is a convex combination of Dirac's deltas
centered in points separated by a distance not smaller than $R$.
\end{proposition}

\begin{pf} Assume by contradiction that there exist
$x^*,y^*\in\operatorname{supp}(\mu_{\infty})$ such that
$|x^*-y^*|<R$. Then, one
can find
suitable neighborhoods $A$ and~$B$ of $x^*$ and $y^*$,
respectively, such that $|x-y|<R$ for all $x\in A$ and $y\in B$.
Hence, $\operatorname{supp}(\theta( \cdot|x,y))\subseteq[\delta
(R), 1]$ for
all $x\in A$  and $y\in B$. Then,
\begin{eqnarray*}
&&\int\!\!\int\!\!\int|x-y|^2\omega\,d\theta(\omega|x,y)\,d
\mu_{\infty}(x)\,d\mu_{\infty}(y)\\
&& \qquad \ge
\delta(R) \int_{ A}\int_{ B} |x-y|^2\,d\mu_{\infty}(x)\,d
\mu_{\infty}(y) > 0 .
\end{eqnarray*}
This clearly contradicts (\ref{Cauchy}).
\end{pf}

It is worth stating the following simple, though important,
consequence of Proposition \ref{propboundedconfidenceasymptotic},
which, in particular, applies to the Gaussian
interaction kernel (\ref{Gaussiandecay}).
\begin{corollary} Suppose that
\[
\bigcup_{\omega_0>0}\bigl\{(x,y)\dvtx \operatorname{supp}\bigl(\theta
(\cdot
|x,y)\subseteq[\omega_0,1]\bigr) \bigr\}=\mathcal X\times\mathcal X .
\]
Then, $\mu_{\infty}=\delta_{x_0}$ where $x_0=\int x\,d\mu_0(x)$.
\end{corollary}

\begin{figure}

\includegraphics{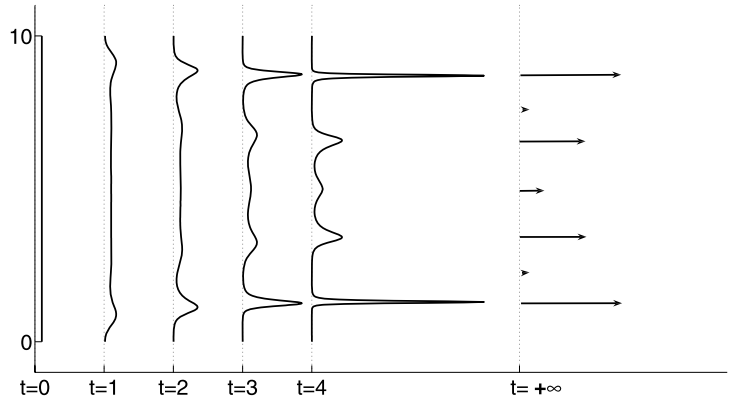}

\caption{Behavior in time of the ODE
solution in $d=1$, with initial condition $\mu_0$ uniform over
$(0,10)$, and $\theta(\cdot|x,y)=\delta_{1/2}\mathbh{1}
_{[0,1]}(|x-y|)+\delta_0\mathbh{1}_{(1,+\infty)}(|x-y|)$. A schematic
representation of the asymptotic distribution ($t=+\infty$), as
presented in \protect\cite{BenNaim}, is reported as well.}\label{figboundedconfidence}
\end{figure}

Figure~\ref{figboundedconfidence} reports numerical simulations of the
mean-field ODE associated to the bounded-confidence model of
Deffuant--Weisbuch, in dimension $d=1$, starting from an initial
condition uniform over the open interval $(0,10)$. Observe that, as
predicted by Proposition \ref{propfinitetimeAC}, the solution remains
absolutely continuous, with bounded density, at any finite time $t$. It
is possible to appreciate the effect of local aggregation forces, which
first lead to the formation of two peaks around the opinion points
$x=1,9$, then of other two smaller peaks around the points $x=3,7$, and
finally of a smaller peak in $x=5$.
As $t$ grows large, the opinion density converges to a convex
combination of Dirac's deltas, as predicted by Proposition \ref
{propboundedconfidenceconvergence}, separated by an inter-cluster
distance of at least $1$, as predicted by Proposition \ref
{propboundedconfidenceasymptotic}. A schematic representation of the
asymptotic opinion distribution, as studied in \cite{BenNaim}, is
plotted as well, presenting some minor clusters between the major ones.
The reader is referred to \cite{BenNaim,Lorenz} for extensive
simulations of this model, and bifurcation studies for the asymptotic
distribution. These results may be interpreted as explaining how
locally aggregating interactions modeling homophily can generate global
fragmentation.

%

We conclude this section by observing that arguments along the lines of
the proofs of Propositions \ref{propboundedconfidenceconvergence} and
\ref{propboundedconfidenceasymptotic}, combined with a standard
martingale convergence theorem, can be used, for every finite
population size $n$, to prove almost sure convergence of the stochastic
system $\mu^n_t$ to a random asymptotic measure $\mu^n_{\infty}$,
consisting of a convex combination of Dirac's deltas separated by a
distance of at least $\sup\{R>0\dvtx  \delta(R)>0\}$.

\section{Concentration around the mean-field dynamics}
\label{sectconcentration} In this section, we finally show that,
as the population size $n$ grows, the stochastic process
$\{\mu^n_t\}$ concentrates around the solution $\{\mu_t\}$ of the
ODE (\ref{ODE}), at an exponential probability rate. Throughout
this section, we shall assume that $\mathcal X\subseteq{\mathbb
{R}}^d$ is
bounded, with $\Delta$ denoting its diameter, and that the
stochastic kernel $\kappa( \cdot| \cdot, \cdot)$ is
globally Lipschitz in the Kantorovich--Wasserstein metric,
that is, that
%
\begin{equation}\label{kappaLip}
\hspace*{5pt}\w(\kappa( \cdot|x,y),\kappa(
\cdot
|x',y') )\le
\frac{L_F}2|(x,y)-(x',y') |  \qquad  \forall x,x',y,y'\in
\mathcal
X
\end{equation}
holds for some finite positive constant $L_F$. Our first step
consists in showing that the operator $F$ inherits the Lipschitz
property from the stochastic kernel
$\kappa( \cdot| \cdot, \cdot)$. The proof of the next
result relies on the duality formula
(\cite{AmbrosioGigliSavare}, (7.1.2))
%
\begin{equation}\label{dualW1}
\w(\mu,\nu)=\sup\{\langle\mu,\varphi\rangle-\langle\nu
,\varphi\rangle\dvtx
\varphi\in\Lip(\mathcal
X) \} ,
\end{equation}
where $\Lip(\mathcal X)$ denotes the set of $1$-Lipschitz
functions on $\mathcal X$.

\begin{lemma}\label{lemmaHLip}
If (\ref{kappaLip}) holds, then
\[
\w(F(\mu),F(\nu))\le L_F\w(\mu,\nu)  \qquad  \forall
\mu,\nu\in\mathcal P_1(\mathcal X) .
\]
\end{lemma}

\begin{pf} First, observe that, for arbitrary function $\varphi\in\Lip
(\mathcal X)$, and $x,y,x'$, $y'  \in\mathcal X$,
\begin{eqnarray*}
\int\varphi(z)\,d\kappa(z|x,y)-\int\varphi
(z)\,d
\kappa(z|x',y')&\le&
\w(\kappa( \cdot|x,y),\kappa( \cdot|x',y') )\\
&\le& \frac{L_F}2|(x,y)-(x',y') |
\\&\le& \frac{L_F}2(|x-x'|+|y-y'| ) ,
\end{eqnarray*}
by (\ref{dualW1})  and (\ref{kappaLip}).
For $\mu,\nu\in\mathcal P(\mathcal X)$, let $\xi\in\mathcal
P_1(\mathcal X\times\mathcal
X)$ be their optimal coupling, that is, the one such that $\int\!\!\int
|x-y|\,d\xi(x,y)=\w(\mu,\nu)$. Then
\begin{eqnarray*}
&&\langle F(\mu),\varphi\rangle - \langle F(\nu),\varphi\rangle\\
&& \qquad  =
\int\!\!\int\!\!\int\!\!\int\!\!\int\varphi(z)\,d\kappa
(z|x,y)\,d\mu(x)\,d\mu(y)\\
&& \qquad  \quad {}- \int\!\!\int\!\!\int\!\!\int\!\!\int\varphi
(z)\,d\kappa(z|x' ,y')\,d\nu(x')\,d\nu(y')\\
&& \qquad  =  \int\!\!\int\!\!\int\!\!\int\!\!\int
\varphi(z) \bigl( d\kappa(z|x,y) - d\kappa(z|x'
,y') \bigr) \,d\xi(x,y)\,d\xi(x' ,y')\\
&& \qquad  \le
\frac{L_F}2\int\!\!\int\!\!\int\!\!\int
(|x-x'|+|y-y'| )\,d\xi(x,y)\,d\xi(x',y')\\
&& \qquad =L_F\w(\mu,\nu) .
\end{eqnarray*}
Hence, the claim follows by applying the duality formula (\ref
{dualW1}) once more.
\end{pf}

Observe that there are three sources of randomness in the system:
the empirical measure of the initial opinions $\mu^n_0$, the
update times $\{T_k\}$, and the agents' interaction. The first two
can be easily dealt with by appealing to the following classical
large deviations results.
\begin{lemma}\label{lemmainitialLD}
For all $\mu_0\in\mathcal
P(\mathcal X)$,
$\eps>0$, it holds that
\[
\limsup_nn^{-1}\log{\mathbb{P}}\bigl(\w(\mu^n_0,\mu_0)\ge\eps
\bigr)\le-\eps^2/2 .
\]
\end{lemma}

\begin{pf}
Sanov's theorem (\cite{ShwartzWeiss}, Theorem~2.14)  and the
Csiszar--Kullback--Pinsker inequality (\cite{Villani}, page~580)  imply that
\begin{eqnarray*}
&&\liminf_n-\frac1n \log{\mathbb{P}}\bigl(\w(\mu
^n_0,\mu
_0 )\ge\eps\bigr) \\
&& \qquad \ge\inf\{H(\nu\parallel\mu_0 )\dvtx  \nu\in\mathcal P(\mathcal
X), \w(\nu,\mu_0)\ge\eps\}\\
&& \qquad \ge\inf\bigl\{\tfrac12\|\nu-\mu_0\|^2\dvtx  \nu\in\mathcal P(\mathcal
X),
\w(\nu,\mu_0)\ge\eps\bigr\}\\
&& \qquad \ge\eps^2/(2\Delta^2)
,
\end{eqnarray*}
where $H(\nu\parallel\mu)$ denoted the relative entropy, and
the last inequality follows from the estimate $W_1(\nu,\mu)\le\Delta
\|\nu-\mu\|$ (\cite{Villani}, Theorem~6.15).
\end{pf}

\begin{lemma}\label{lemmaPoissonLD}
For $t\in{\mathbb{R}}^+$, let $\varsigma(t):=\sup\{k\in{\mathbb
{Z}}^+\dvtx  T_k\le t\}$.
For all $\tau\in{\mathbb{R}}^+$, $a\ge1$ it holds
\begin{eqnarray*}
&\displaystyle\limsup_nn^{-1}\log{\mathbb{P}}\bigl(\sup\bigl\{t-T_{\varsigma(t)}\dvtx
0\le
t\le\tau\bigr\}\ge\eps\bigr)\le-\eps^2/\tau,&
\\
&\displaystyle\limsup_nn^{-1}\log{\mathbb{P}}\bigl(\varsigma(\tau)\ge a\tau
n\bigr)\le
-(a-1)^2\tau.&
\end{eqnarray*}
\end{lemma}

\begin{pf} The first statement follows, for example, from \cite{ShwartzWeiss}, Theorem~5.1.
The second one, for example, from \cite{ShwartzWeiss}, Example~1.13.
\end{pf}

We are now left with the third source of randomness, originated by the
selection of the interacting agents, and their actual interaction.
Observe that, in the right-hand side of the duality formula (\ref
{dualW1}), one may restrict the supremization
to the test functions $\varphi$ belonging to
\[
\Lipd:=\{\varphi\in\Lip(\mathcal Y)\dvtx  |\varphi(x)|\le\Delta
/2 \} ,
\]
where $\mathcal Y$ is a hypercube of edge-length $\Delta$ containing
$\mathcal
X$, and $\mu,\nu$ are naturally identified as elements of $\mathcal
P(\mathcal Y)$.
The following result shows that the set $\Lipd$ can be approximated in
the infinity norm by not-too-large a set of functions.
\begin{lemma}\label{lemmaLipschitzapproximation}
Let $\mathcal X\subseteq{\mathbb{R}}^d$ be compact and convex. Then,
for all $\eps
\in\,]0,\Delta/2]$, there exists a finite set $\mathcal H_{\eps
}\subseteq
\Lipd$ such that
$|\mathcal H_{\eps}|\le\frac{2\sqrt d+1}{6}\frac{\Delta}{\eps}3^{
(\fraca{\Delta}{\eps}(\sqrt d+1) )^d}$,
and
\[
\min\{\|h-\varphi\|\dvtx  h\in\mathcal H_{\eps} \}\le\eps
 \qquad
\forall\varphi\in\Lipd.
\]
\end{lemma}

\begin{pf} With no loss of generality, we shall restrict to the case
$\mathcal X\subseteq\mathcal Y=[0,\Delta]^d$. We introduce a discretization
operator $\Phi\dvtx \Lipd\to\Lipd$ as follows. Let $\eta:=\eps/(\sqrt
d+1/2)$ and define $\mathcal
J:=\{0,1,\ldots,\lfloor\Delta/\eta\rfloor\}$. For any
$\varphi\in\Lipd$  and $\mb j\in\mathcal J^d$, let $k(\mb j)=i\in
\mathcal J$
iff $\varphi(\mb j\eta)\in[-1/2+\eta i,-1/2+\eta(i+1)[$. Observe
that, since $\varphi$ is $1$-Lipschitz, one has
%
\begin{equation}\label{kj}
\sum_{1\le l\le d}|j_l-j'_l|\le
1  \quad \Longrightarrow \quad |k(\mb j)-k(\mb j')|\le1 .
\end{equation}
Then,
define $\Phi(\varphi)=h$, by putting, for all $x\in\prod_{1\le
l\le d}[j_l\eta,(j_l+1)\eta]$,
\[
h(x)=\prod_{1\le l\le d}\biggl(\bigl(k(\mb j+\delta_l)-k(\mb j)
\bigr)(x_l-j_l\eta)+\eta k(\mb j)-\frac12+\frac{\eta}2 \biggr) .
\]

Thanks to (\ref{kj}), one has that $\Phi(\varphi)\in\Lipd$ for all
$\varphi\in\Lipd$.
Moreover, for all $\mb j\in\mathcal J^d$, one has $|\Phi(\varphi
)(\mb
j\eta)-\varphi(\mb j\eta)|\le\frac\eta2$.
Observe that, for all $x\in[0,\Delta]^d$, there exists $\mb j(x)\in
\mathcal J^d$ such that $|x-\eta\mb j|\le\sqrt d\eta/2$. Therefore,
\begin{eqnarray*}
|\Phi(\varphi)(x)-\varphi(x)|&\le&
|\Phi(\varphi)(\mb j\eta)-\varphi(\mb j\eta)|+|\Phi(\varphi)(\mb
j\eta)-\Phi(\varphi)(x)|\\
&&{}+|\varphi(\mb j\eta)-\varphi(x)|\\
&\le&
\eta/2+2|\mb j\eta-x|\\
&\le&\eta\bigl(\sqrt d+1/2\bigr) ,
\end{eqnarray*}
so that the second part of the claim follows by substituting the value
of $\eta$.

It remains to estimate the cardinality of $\mathcal H_{\eps}:=\Phi
(\Lipd)$.
To see that, first observe that $k(\mb0)$ can take at most $\Delta
/\eta$ values. On the other hand, it follows from (\ref{kj}) that,
given $k(\mb j)$, $k(\mb j+\delta_l)$ can assume at most three
different values, for all $1\le l\le d$. This implies that
\begin{eqnarray*}
|\mathcal H_{\eps}|\le\frac{\Delta}{\eta
}3^{(\Delta/\eta
+1)^d-1}&=&
\frac{\Delta}{\eps}\frac{2\sqrt d+1}{6}3^{((\sqrt
d+1/2)\Delta/\eps+1 )^d}\\&\le&
\frac{\Delta}{\eps}\frac{2\sqrt d+1}{6}3^{((\sqrt d+1)\Delta
/\eps)^d} ,
\end{eqnarray*}
the last inequality following since $1\le\Delta/(2\eps)$.
\end{pf}

We can now estimate the error incurred when using an Euler
approximation of some future value of the empirical density process,
centered on its current value.

\begin{lemma}\label{lemmaAzuma}
For $k\in{\mathbb{Z}}^+$, $n\in{\mathbb{N}}$  and $\sigma\in[0,1]$,
\[
{\mathbb{P}}\bigl(\w\bigl(\ov\sigma M_{k}+\sigma F(M_{k}),M_{k+\lfloor\sigma
n\rfloor} \bigr)\ge K\Delta\sigma^2 \bigr)\le\rho,
\]
where $\ov\sigma=1-\sigma$, $K=K_F+1$, with $K_F$ being the
Lipschitz constant of $F$ on $\mathcal P(\mathcal X)$ in the variational
distance, and
%
\begin{equation}\label{rhodef}\rho:=\frac{4\sqrt d+2}{K\sigma
^2}\exp\biggl(
\biggl(\frac{12}{K\sigma^2}\bigl(\sqrt d+1\bigr) \biggr)^d\log3-\frac{K^2\sigma
^3}{2^7}n \biggr) .
\end{equation}
\end{lemma}

\begin{pf}
First, observe that the following control of the increments holds:
%
\begin{equation}\label{TVincrement}
\|M_{k+1}-M_{k}\|\le2/n .
\end{equation}
Define $w:=\lfloor\sigma n\rfloor$  and $\eps:=K\Delta\sigma^2$.
Also, for $\varphi\in\Lipd$, define
\[
Z^{(\varphi)}_j:=\langle M_{k+j}-M_{k},\varphi\rangle-\frac1n\sum
_{0\le
i< j}\langle F(M_{k+i} )-M_{k+i},\varphi\rangle
\]
for $j=0,\ldots,w$, and
\[
V^{(\varphi)}:=\biggl\langle M_{k+w}-\biggl(1-\frac wn\biggr)M_{k}-\frac
wnF(M_{k}),\varphi
\biggr\rangle-Z^{(\varphi)}_w .
\]
It follows from (\ref{TVincrement}) that $\|M_{k+j}-M_{k}\|\le2j/n$.
Hence, 
%
\begin{eqnarray}\label{Vphiestimate}
\bigl|V^{(\varphi)} \bigr|
&=&n^{-1} \biggl|\sum_{0\le j<w}\langle F(M_{k+j})-F(M_{k}),\varphi
\rangle
-\sum_{0\le j<w}\langle M_{k+j}-M_{k},\varphi\rangle\biggr|
\nonumber\\
&\le&
n^{-1}\sum_{0\le j<w}
\bigl(\|F(M_{k+j})-F(M_{k})\|+\|M_{k+j}-M_{k}\| \bigr)\|\varphi\|\nonumber
\\[-8pt]
\\[-8pt]
&\le&n^{-1}\sum_{0\le j<w}K\frac{2j}n\|\varphi\|\nonumber\\
&\le&\eps/2 ,
\nonumber
\end{eqnarray}
the last inequality following from the fact that $\|\varphi\|\le
\Delta/2$.
Observe that, for all
$\varphi\in\Lip(\mathcal X)$, $Z^{(\varphi)}_0=0$, while $\{
Z^{(\varphi
)}_j\dvtx  0\le j\le w\}$ is a martingale.
Moreover, (\ref{TVincrement}) provides the following control on the increments:
%
\begin{eqnarray}\label{increments}
 \bigl|Z^{(\varphi)}_{j+1}-Z^{(\varphi)}_j\bigr|&\le&
|\langle M_{k+j+1}-M_{k+j},\varphi\rangle| +n^{-1}|\langle
F(M_{k+j})-M_{k+j},\varphi\rangle|\nonumber\\
&\le&\|M_{k+j+1}-M_{k+j}\|\|\varphi
\|+n^{-1}\|F(M_{k+j})-M_{k+j}\|\|\varphi\|
\\
&\le&4n^{-1}\|\varphi\|
.\nonumber
\end{eqnarray}

Let $\mathcal H:=\mathcal H_{\eps/12}\subseteq\Lip(\mathcal X)$ be
as in Lemma \ref
{lemmaLipschitzapproximation}.
By first applying the union bound, and then the Hoeffding--Azuma
inequality (\cite{AlonSpencer}, Theorem~7.2.1)  the probability of the event
%
$E:=\bigcup_{h\in\mathcal H}\{|Z^{(h)}_w|\ge\eps/4 \}$
can be estimated as follows:
%
\begin{equation}{\mathbb{P}}(E )\le|\mathcal H|{\mathbb
{P}}\bigl(\bigl|Z^{(h)}_w\bigr|\ge\eps/4 \bigr)
\le2|\mathcal H|\exp\biggl(-\frac{\eps^2n^2}{2^7w\Delta^2} \biggr)
.\label{P(E)estimate}
\end{equation}
%
Now, Lemma \ref{lemmaLipschitzapproximation} and (\ref{increments})
imply that
\[
Z^{(\varphi-h)}_w\le3\frac wn\|\varphi-h\|\le3\sigma\frac{\eps
}{12}\le\frac{\eps}4
\]
for some $h\in\mathcal H_{\eps/12}$. Hence, if $E$ does not occur, then
%
\begin{equation}\label{Zwest}
\bigl|Z_w^{(\varphi)}\bigr|\le\min\bigl\{
\bigl|Z_w^{(h)}\bigr|+\bigl|Z_w^{(\varphi-h)}\bigr|\dvtx  h\in\mathcal H \bigr\}\le\frac
\eps2
\end{equation}
for every $\varphi\in\Lipd$.
By combining (\ref{Vphiestimate}), (\ref{P(E)estimate})  and (\ref
{Zwest}), one gets
\begin{eqnarray*}
{\mathbb{P}}\bigl(\w\bigl(M_{k+w},\ov\sigma M_{k}+\sigma F(M_{k}) \bigr)\ge
\eps
\bigr)
&=&{\mathbb{P}}\bigl(\sup\bigl\{Z^{(\varphi)}_w+V^{(\varphi)}\bigr\}\ge\eps
\bigr)\\
&\le&{\mathbb{P}}\biggl(\sup\bigl\{Z^{(\varphi)}_w\bigr\}\ge\frac34\eps\biggr)\\
&\le&2|\mathcal H|\exp\biggl(-\frac{\eps^2n^2}{2^7w\Delta^2} \biggr) ,
\end{eqnarray*}
and the claim follows upon substituting the expressions for $w$ and
$\eps$, and applying Lemma \ref{lemmaLipschitzapproximation}.
\end{pf}

We are now ready to prove point (b) of Theorem \ref{theomain}.
Let $L:=L_F-1$  and $K=K_F+1$, where $L_F$ and $K_F$ are the global
Lipschitz constants of $F$ on~$\mathcal P(\mathcal X)$ in the
Kantorovich--Wasserstein distance,
and in the variational distance, respectively.
Let us fix some $\eps>0$, $\tau>0$, and introduce the quantities
\[
\sigma:=\frac{L\eps}{2\Delta L+3K\Delta e^{2L\tau}} , \qquad
w=\lfloor\sigma/n\rfloor.
\]
Without any loss of generality, let us assume that $\sigma\in\,]0,1]$,
and put $\ov\sigma=1-\sigma$. Further, let $\rho$ be as in
(\ref{rhodef}), and define
%
\begin{equation}\label{alphadef}
\alpha_0=e^{-2L\tau}\eps/2, \qquad
\alpha_{i+1}= (1+\sigma
L )\alpha_i+\tfrac32K\Delta\sigma^2 , \qquad  i\in{\mathbb
{Z}}^{+} .
\end{equation}
Solving the iterative equation above, one obtains the estimate
%
\begin{equation}\label{alphaestimate}
\alpha_i=(1+\sigma
L )^{i}\biggl(\alpha_0+\frac{3K\Delta\sigma}{2L} \biggr)-\frac
{3K\Delta\sigma}{2L}\le
e^{\sigma Li}\biggl(\alpha_0+\frac{3K\Delta\sigma}{2L} \biggr) .
\end{equation}
For
$i\in{\mathbb{Z}}^+$, consider the random variable
$\Gamma^n_i:=\w(M_{iw},\mu_{\sigma i} )$, and the events
$A_i:=\{\Gamma^n_i\ge\alpha_i\}$, $B_i:=\bigcup_{0\le j\le i}A_j$.
We shall prove by induction that
%
\begin{equation}\label{induction}
\mathbb{P}(B_i )\le(i+1)\rho
\end{equation}
for all $i\in{\mathbb{Z}}_+$. First,
it follows from Lemma \ref{lemmainitialLD} that (\ref{induction})
holds with $i=0$, for sufficiently small $\eps$  and sufficiently
large $n$. Then, for any nonnegative integer $i$, consider the
intermediate measures
\[
\lambda:=\ov\sigma M_{wi}+\sigma F(M_{wi} ) , \qquad  \nu
:=\ov\sigma\mu_{\sigma i}+\sigma F(\mu_{\sigma i} ) .
\]
From the duality formula (\ref{dualW1}), and Lemma \ref{lemmaHLip},
one has
%
\begin{equation}\label{triangle2}
\w(\lambda,\nu)\le(\ov\sigma+\sigma L_F) \Gamma
^n_{i}=(1+\sigma L) \Gamma^n_{i} .
\end{equation}
Furthermore, since $\{\mu_t\}$ is a solution of the ODE (\ref
{generalODE}), it follows from (\ref{Phitcont}), and the estimate $\w
(\mu,\nu)\le\Delta/2\|\mu-\nu\|$,
%
\begin{equation}\label{W1estimate}
\|\mu_t-\mu_{\sigma i}\|\le
2(t-s )
, \qquad  \w(\mu_t,\mu_s)\le\Delta(t-s)
\end{equation}
for all $t\ge s$.
From the duality formula (\ref{dualW1}), the fact that $\{\mu_t\}$
solves the ODE (\ref{generalODE}), and (\ref{W1estimate}), one gets
the estimate
%
\begin{eqnarray}\label{triangle3}
\w\bigl(\nu,\mu_{\sigma(i+1)} \bigr) & = &
\sup\bigl\{\bigl\langle\mu_{\sigma(i+1)},\varphi\bigr\rangle-\langle\nu
,\varphi\rangle
\dvtx  \varphi\in\Lipd\bigr\}\nonumber\\[2pt]
&\le&
\int_{\sigma i}^{\sigma(i+1)} \sup\{\langle
F(\mu_t)-\mu_t-F(\mu_{\sigma i})+\mu_{\sigma i},\varphi\rangle\dvtx
\varphi\in\Lipd
\}\,dt\nonumber\\[2pt]
&\le&
\frac{\Delta}2K\int_{\sigma i}^{\sigma(i+1)}\|\mu_t-\mu
_{\sigma i}\|\,dt
\\[2pt]
&\le&
\Delta K\int_{\sigma i}^{\sigma(i+1)}(t-\sigma i )\,dt
\nonumber\\
&=&\Delta K\sigma^2/2 .\nonumber
\end{eqnarray}
From the triangle inequality, (\ref{triangle2})  and (\ref
{triangle3}), one finds that
%
\begin{eqnarray}\label{W1triangle}
 \Gamma^n_{i+1}&\le&\w\bigl(M_{w(i+1)},\lambda\bigr)+\w
(\lambda,\nu
)+\w\bigl(\nu,\mu_{\sigma(i+1)} \bigr)\nonumber
\\[-6pt]
\\[-6pt]
&\le&\w\bigl(M_{w(i+1)},\lambda\bigr)+K\Delta\sigma^2/2+(1+\sigma L
)\Gamma^n_{i} .
\nonumber
\end{eqnarray}
Therefore, (\ref{W1triangle}), the inductive hypothesis (\ref
{induction}), (\ref{alphadef})  and Lemma \ref{lemmaAzuma}, imply that
\begin{eqnarray*}
{\mathbb{P}}(B_{i+1})&=&
{\mathbb{P}}(B_i^c\cap A_{i+1} )+{\mathbb{P}}(B_i)\\
&\le&
{\mathbb{P}}\bigl(B_i^c\cap\bigl\{\w\bigl(M_{w(i+1)},\lambda\bigr)>\alpha
_{i+1}-K\Delta
\sigma^2/2-(1+\sigma L )\alpha_{i} \bigr\} \bigr)\\
&&{}+(i+1)\rho\\
&\le&{\mathbb{P}}\bigl(\w\bigl(M_{w(i+1)},\lambda\bigr)> K\Delta\sigma
^2\bigr)+(i+1)\rho\\
&\le&(i+2)\rho.
\end{eqnarray*}
Hence, (\ref{induction}) holds for all $i\in{\mathbb{Z}}^+$.

Observe that, if $iw-w/2\le k\le iw+w/2$, then
%
\begin{eqnarray}\label{Mkmukn}
 \w(M_{k},\mu_{k/n})&\le&\w(M_{wi},\mu
_{\sigma i})+\w(M_{wi},M_{k})+\w(\mu_{k/n},\mu_{\sigma i})\nonumber
\\[-8pt]
\\[-8pt]
&\le&\Gamma^n_i+\Delta\sigma.
\nonumber
\end{eqnarray}
Now, recall the definition of $\varsigma(t)$ given in Lemma \ref
{lemmaPoissonLD}, and consider the events $C:=\{\varsigma(\tau)\le
\frac32n\tau\}$  and $D:=\{\sup\{|t-T_{\varsigma(t)}|\dvtx  t\in
[0,\tau]\}\le\eps/(4\Delta) \}$. Observe that $C$ implies
that, for all $t\le\tau$,
%
\begin{equation}\label{boring}
\iota(t):=\biggl\lfloor\frac{\varsigma
(t)}{\lfloor
\sigma n\rfloor}+\frac12 \biggr\rfloor\le\frac{3\tau n/2}{\sigma
n-1}+\frac12\le\frac{2\tau}{\sigma} .
\end{equation}
It follows from (\ref{Mkmukn}), (\ref{W1estimate}), (\ref
{alphaestimate})  and (\ref{boring}), that, if the event $B^c_{\lfloor
2\tau\sigma\rfloor}\cap D\cap C$ occurs, then, for all $t\in[0,\tau
]$, the following estimate holds:
\begin{eqnarray*}
\w(\mu^n_t,\mu_t )&=&\w\bigl(M_{\varsigma
(t)},\mu
_t \bigr)\\
&\le&\w\bigl(M_{\varsigma(t)},\mu_{T_{\varsigma(t)}} \bigr)+\w
\bigl(\mu_{T_{\varsigma(t)}},\mu_t \bigr)\\
&\le&\Gamma^n_{\iota(t)}+\Delta\sigma+\Delta\bigl|t-T_{\varsigma
(t)}\bigr|\\
&\le&\alpha_{\iota(t)}+\Delta\sigma+\eps/4\\
&\le&e^{\sigma L\iota(t)}+\Delta\sigma+\eps/4\\
&\le&e^{2L\tau}\biggl(\alpha_0+\frac{3K\Delta\sigma}{2L}
\biggr)+\Delta\sigma+\eps/4\\
&=&\eps,
\end{eqnarray*}
where the last equality follows by substituting the expressions for
$\sigma$ and $\alpha_0$.
For sufficiently small $\eps$  and large $n$, Lemma \ref
{lemmaPoissonLD} implies that ${\mathbb{P}}(C\cap D)\ge1-\rho$.
Therefore, using (\ref{induction}), one gets that
\[
{\mathbb{P}}\bigl(\sup\{\w(\mu^n_t,\mu_t )\dvtx  t\in[0,\tau]\}
>\eps
\bigr)\le{\mathbb{P}}\bigl(B_{\iota(\tau)}\bigr)+{\mathbb{P}}(C^c\cup
D^c)\le(2\tau/\sigma
+2)\rho,
\]
from which point (b) of Theorem \ref{theomain} follows.

%

\printaddresses

\end{document}